\documentclass[11pt,reqno]{amsart}

\usepackage[utf8]{inputenc}
\usepackage[english]{babel}
\usepackage{amsmath,amssymb,amsthm}
\usepackage{graphicx}
\usepackage{tikz}
\usepackage{caption}
\usepackage{subcaption}
\usepackage{booktabs}
\usepackage{hyperref}

\theoremstyle{plain}
\newtheorem{theorem}{Theorem}[section]
\newtheorem{proposition}[theorem]{Proposition}
\newtheorem{lemma}[theorem]{Lemma}

\newtheorem{conjecture}[theorem]{Conjecture}

\theoremstyle{definition}
\newtheorem{definition}[theorem]{Definition}
\newtheorem{observation}[theorem]{Observation}

\theoremstyle{remark}
\newtheorem{remark}[theorem]{Remark}

\begin{document}

\title{Maps on Surfaces and the Tabulation of Knots and Links
in the Thickened Torus Through Ten Crossings}

\title[Genus-One Tabulation Through Ten Crossings]{Maps on Surfaces and the Tabulation of Knots and Links
in the Thickened Torus Through Ten Crossings}

\author{Alexander Omelchenko}
\address{Constructor University Bremen, Campus Ring 1, 28759 Bremen, Germany}
\email{aomelchenko@constructor.university}

\maketitle

\begin{abstract}
We give a combinatorial tabulation of knots and links in the thickened torus
$T^2\times I$ based on the theory of maps on surfaces: cellular $4$--regular
torus projections are encoded by permutation pairs, and unsensed equivalence
classes are enumerated completely and without duplication by canonical
representatives.  The method is validated against published genus-one tables
for $N\le 5$ and extended to $N=6,7,8,9,10$, producing, to our knowledge, the first complete tabulation
of prime knot and link diagram types in the fixed thickened torus $T^2\times I$
at these crossing numbers, under the conventions stated below --- more than
$1.3$ million knot and link types in total. The computation yields three proved structural results on
straight-ahead components, on bigon faces, and on the parity of the $a$-span
of the genus-one bracket, together with several conjectures, including a $4N$
bound on the $a$-span in the spirit of the Kauffman--Murasugi--Thistlethwaite
theorem.  A reference implementation and machine-readable datasets are provided.
\end{abstract}


\section*{Introduction}\label{sec:introduction}

Tabulation has long played an important role in knot theory, both as a source of examples and as a guide for further developments.
In the classical case this tradition goes back to Tait and, in modern form, leads to very large knot tables; see, for example,
\cite{Tait1877, HosteThistlethwaiteWeeks1998}.

Virtual knots and, more generally, knots and links in thickened surfaces $\Sigma_g\times I$ provide a natural extension of this
program beyond $S^3$.
Virtual knot theory was introduced by Kauffman~\cite{Kauffman1999}.
A standard topological reformulation identifies virtual links with stable equivalence classes of link diagrams on oriented
surfaces, or equivalently with links in thickened surfaces modulo stabilization and destabilization~\cite{KamadaKamada2000, CarterKamadaSaito2002}.
Kuperberg's theorem implies that every virtual link admits a minimal-genus representative and that the supporting surface of such
a representative is unique up to homeomorphism~\cite{Kuperberg2003}.
In particular, genus-one virtual knots and links may be studied through diagrams on the thickened torus $T^2\times I$.

For the torus, explicit low-crossing classifications are already available.
The first systematic genus-one lists at very small complexity, together with surface-sensitive bracket-type invariants, were developed
in~\cite{Omelchenko2007, Grishanov2009a, Grishanov2009b, Morton2009}.
Later, Akimova and Matveev classified knots in $T^2\times I$ admitting diagrams with at most four crossings~\cite{AkimovaMatveev2012};
Akimova and Matveev produced a genus-one virtual-knot table up to five classical crossings~\cite{AkimovaMatveev2014};
and Akimova, Matveev and Tarkaev extended the genus-one program further, including a complete table of prime links in $T^2\times I$
with crossing number five~\cite{AkimovaMatveevTarkaev2018, AkimovaMatveevTarkaev2020}.
Complementary tabulations in the broader virtual setting include Green's table of virtual knots up to six crossings~\cite{GreenTableVirtualKnots},
Chrisman's classification of slice status up to five crossings~\cite{Chrisman2022Milnor}, and the tabulation of almost classical knots
by Boden et al.~\cite{BodenGaudreauHarperNicasWhite2017}.
These works provide the natural benchmark and context for the present paper.

A persistent difficulty in surface tabulation lies in the projection stage.
A common approach starts from abstract $4$--regular graphs, then classifies their inequivalent embeddings into the surface up to
homeomorphism, and only afterwards assigns crossing information to obtain diagrams.
While conceptually straightforward, this becomes difficult to scale: the classification of embeddings quickly proliferates into cases,
and knot and link projections are naturally separated because the component structure becomes visible only after the embedding has been fixed.

Our starting point is to treat the embedding itself as the primary combinatorial object,
following the classical theory of maps on surfaces.
A cellular embedding of a connected graph into a closed orientable surface is a map, and a cellular $4$--regular torus projection is
therefore a $4$--regular map on $T^2$.

The dart-permutation encoding used below is classical: an involution $\alpha$ pairs the two darts of each edge, while a
permutation $\sigma$ records the cyclic order of darts around each vertex.
This viewpoint goes back at least to Walsh--Lehman and is standard in the theory of maps on surfaces; see also
\cite{WalshLehman1972,LandoZvonkin2004}.
A closely related permutation-orbit framework was used by Coquereaux and Zuber to enumerate immersed circles in orientable surfaces,
with one-componentness and genus imposed as cycle conditions on permutations~\cite{CoquereauxZuber2016}.

Our use of this framework is different in scope: we use it as the projection stage of a tabulation pipeline for prime knots and links
in the fixed thickened torus $T^2\times I$.
In this language, equivalence up to homeomorphism becomes a finite conjugacy relation on permutation pairs, and several projection-level
structures relevant to knot theory become explicit: the number of straight-ahead components can be read off directly, and monogon faces
corresponding to immediate Reidemeister~I reductions are detected combinatorially.
Thus the quotient by surface homeomorphisms is performed already at the projection stage, before crossing data and invariants are introduced.

The present work builds on earlier results on quartic maps on surfaces, including one-face quartic maps, toroidal map families, and rooted
quartic maps without planar monogons~\cite{Krasko_Omelch_4_reg_one_face_maps, KraskoOmelchenko2019PartI, KraskoOmelchenko2019PartII, KraskoOmelchenko2021NoPlanarLoops}.
Here these map-theoretic tools are used not merely to count maps, but to construct a practical and reproducible classification pipeline
for knots and links in the thickened torus.

More precisely, for each crossing number $N$ we enumerate loopless cellular $4$--regular torus projections without monogons, quotient them
by unsensed equivalence via canonical representatives, and apply explicit projection-level primeness witnesses.
For each fixed projection we encode crossing choices by bit assignments, formulate a local Reidemeister~II reduction criterion along a
bigon face directly in terms of these bits, and compute the genus-one generalized bracket entirely within the same permutation model,
without drawing diagrams in a fundamental polygon.
The resulting pipeline is validated against the published genus-one classifications for crossing numbers $N\le 5$ under explicitly stated
comparison conventions.

We then extend the computation to $N=6,7,8,9,10$, obtaining, to our knowledge, the first complete genus-one tabulation in this range.
The released datasets contain unsensed canonical projection classes, diagram representatives, and the corresponding invariants, all generated
by a public reference implementation~\cite{OmelchenkoTorusMaps}.
Through $N=10$ the tabulation contains $1\,357\,206$ knot and link types.

Beyond the tables themselves, the computation leads to several structural results and conjectures.
We prove, by an explicit lattice construction, that the maximum number of straight-ahead components of a prime torus projection is at least
$\lfloor N/2\rfloor+2$ for even $N$, and the detailed projection data analyzed below support the equality
$c_{\max}(N)=\lfloor N/2\rfloor+2$ throughout the range $N\le 10$.
We also prove that the $a$-span of the genus-one bracket is always even, and that a $4$--regular torus projection with $N$ vertices has at most
$N-1$ bigon faces; the latter bound is attained for every $N\le 10$ in the detailed projection data analyzed below.
The computations further indicate, in the analyzed invariant data through $N\le 10$, that the $a$-span satisfies the upper bound $4N$, in analogy with the classical
Kauffman--Murasugi--Thistlethwaite theorem, and that every computed equi-quadrilateral prime projection is a link projection.
These findings, together with the resulting conjectures, are collected in Section~\ref{sec:observations}.

The paper is organized as follows.
Section~\ref{sec:maps} recalls the permutation model for maps and the projection-level structures needed later.
Section~\ref{sec:projections} specializes to the torus and proves that canonicalization yields a complete and duplication-free enumeration of
unsensed prime projection classes, together with benchmark counts up to $N=10$.
Section~\ref{sec:diagrams} passes from projections to diagrams, establishes the bigon reduction criterion, and formulates the genus-one bracket
invariants entirely in the permutation model; the implementation is validated there against the published classifications for $N\le 5$.
Section~\ref{sec:observations} records structural consequences and conjectures suggested by the computed data.
Finally, Section~\ref{sec:outlook} discusses possible extensions beyond genus one.
Algorithmic details and reproducibility notes are collected in Appendices~\ref{app:proj_algorithms} and~\ref{app:diagram_algorithms}.

\section{Maps on surfaces}\label{sec:maps}

Throughout the paper we work with \emph{closed, connected, orientable} surfaces.
For $g\ge 0$ we write $\Sigma_g$ for the closed orientable surface of genus~$g$ (so $\Sigma_1=T^2$).
At the level of maps we allow finite connected multigraphs, including multiple edges and, unless explicitly excluded,
loop edges. The treatment of loop edges will be made explicit, since for the genus-one projection lists considered later
we impose a looplessness constraint.

\begin{definition}[Map on a surface]\label{def:map_on_surface}
A \emph{(topological) map} on $\Sigma_g$ is a cellular embedding of a connected multigraph $G$ into $\Sigma_g$,
i.e.\ an embedding such that each connected component of $\Sigma_g\setminus G$ is homeomorphic to an open disc.
Equivalently, the embedding induces a CW decomposition of $\Sigma_g$ whose $0$-, $1$-, and $2$-cells are the vertices,
edges, and complementary regions of the embedded graph.
\end{definition}

In the knot-theoretic language of thickened surfaces, a \emph{projection} on $\Sigma_g$ is an embedded $4$--regular graph
without over/under data.
In this paper, an \emph{essential projection} means a \emph{cellular} projection in the sense of
Definition~\ref{def:map_on_surface}.
A \emph{diagram} is a projection together with crossing information at each vertex.

We recall the classical dart-permutation encoding of orientable maps.
This material is included to make the later algorithms self-contained for readers coming from knot theory rather than map enumeration.
The formalism goes back at least to Walsh--Lehman and is standard in the theory of maps on surfaces;
see \cite{WalshLehman1972} and \cite[\S 1.3]{LandoZvonkin2004}.
A closely related use of the same language for immersed curves appears in \cite{CoquereauxZuber2016}.
Let $H$ be the dart set, so $|H|=2|E|$.

\smallskip
\noindent\textbf{Convention (right action).}
Permutations act on the right: for permutations $p,q$ we write $pq$ for the composition “first $p$, then $q$”,
so $(pq)(h)=q(p(h))$.

\begin{definition}[Permutation encoding]\label{def:perm_encoding}
A \emph{labelled orientable map} is encoded by a pair of permutations $(\alpha,\sigma)$ on $H$ such that:
\begin{enumerate}
\item $\alpha$ is a fixed-point-free involution (it pairs the two darts of each edge);
\item the cycles of $\sigma$ record the cyclic order of darts around vertices;
\item the subgroup $\langle \alpha,\sigma\rangle$ acts transitively on $H$ (this is connectedness).
\end{enumerate}
The associated \emph{face permutation} is
\[
\varphi:=\sigma\alpha .
\]
\end{definition}
This permutation encodes the oriented boundary walks of faces.

In this standard encoding, the cycles of $\sigma$ are the vertices, the $2$--cycles of $\alpha$ are the edges, and the
cycles of
\[
\varphi=\sigma\alpha
\]
are the oriented face-boundary walks. Thus
\begin{equation}\label{eq:classical_map_counts}
|V|=\#\mathrm{cyc}(\sigma),\qquad
|E|=\#\mathrm{cyc}(\alpha)=\frac{|H|}{2},\qquad
|F|=\#\mathrm{cyc}(\varphi).
\end{equation}

For clarity, let us spell out the face step with our right-action convention.  Starting from a dart $h$, one first turns
to the next dart at the same vertex and then crosses the corresponding edge, so one boundary step is
$h\mapsto(\sigma\alpha)(h)=\alpha(\sigma(h))$.  Iterating this step traces exactly one oriented face boundary.

Similarly, the underlying embedded graph is connected if and only if $\langle\alpha,\sigma\rangle$ acts transitively on $H$:
applying $\alpha$ moves across an edge, while applying $\sigma$ moves between darts incident to the same vertex.
These are standard facts of the permutation model for maps; we use them below as conventions of the encoding.
For a cellular embedding on $\Sigma_g$ we therefore recover the genus from Euler's formula:
\[
\#\mathrm{cyc}(\sigma)-\#\mathrm{cyc}(\alpha)+\#\mathrm{cyc}(\sigma\alpha)=2-2g.
\]

\begin{figure}[ht]
\centering
\includegraphics[scale=1.3]{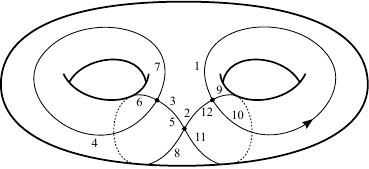}
\caption{A one-face $4$-regular map (genus $2$ example used in the text).
Vertices are shown as crossings; darts are labelled $1,\dots,12$.}
\label{fig:intertwined}
\end{figure}

\medskip
\noindent\textbf{A worked example.}
Figure~\ref{fig:intertwined} shows a one-face $4$-regular map on an orientable surface of
genus~$2$, encoded by permutations on $H=\{1,\dots,12\}$:
\[
\alpha=(1,10)(2,12)(3,5)(4,7)(6,8)(9,11),\qquad
\sigma=(1,12,10,9)(2,5,8,11)(3,7,6,4).
\]
Then $\varphi=\sigma\alpha$ is the single cycle
$(1,2,3,4,5,6,7,8,9,10,11,12)$, hence $F=1$. Since $E=6$ and $V=3$, we obtain
$V-E+F = 3-6+1=-2$, hence $g=2$.
For the torus projections that occupy the rest of the paper the same
encoding applies with the genus condition $V-E+F=0$.

\medskip
Changing the dart labels means applying a bijection $\pi:H\to H$, i.e.\ a permutation $\pi\in S_H$. At the level of encodings this is simultaneous conjugation:
\[
(\alpha,\sigma)\longmapsto (\pi\alpha\pi^{-1},\,\pi\sigma\pi^{-1}).
\]
This is the intrinsic way of “forgetting labels”.

Beyond relabelling, one further identifies maps that differ by a homeomorphism of the supporting surface.
Two conventions are standard: \emph{sensed} equivalence allows only orientation-preserving
homeomorphisms, while \emph{unsensed} equivalence allows orientation reversal as well.
Since reversing the orientation of the surface reverses the cyclic order of darts at each vertex,
at the level of encodings it replaces $\sigma$ with $\sigma^{-1}$ while leaving $\alpha$ unchanged.

\begin{definition}[Sensed and unsensed equivalence]\label{prop:maps_conjugacy_equivalence}
Let $(\alpha_1,\sigma_1)$ and $(\alpha_2,\sigma_2)$ be labelled encodings of orientable maps.
They define the same \emph{sensed} unlabelled map if they are related by simultaneous relabelling,
\[
(\alpha_2,\sigma_2)=(\pi\alpha_1\pi^{-1},\pi\sigma_1\pi^{-1})
\quad\text{for some }\pi\in S_H.
\]
They define the same \emph{unsensed} unlabelled map if, in addition, we also allow reversal of the orientation of the supporting
surface, which replaces the vertex rotation by its inverse:
\[
(\alpha_2,\sigma_2)=(\pi\alpha_1\pi^{-1},\pi\sigma_1\pi^{-1})
\quad\text{or}\quad
(\alpha_2,\sigma_2)=(\pi\alpha_1\pi^{-1},\pi\sigma_1^{-1}\pi^{-1}).
\]
\end{definition}

This is the standard relabelling action on labelled permutation encodings of maps.
The same orbit language is used, for example, in the enumeration of immersed curves by Coquereaux--Zuber~\cite{CoquereauxZuber2016}.

We now record two projection-level structures that will be used throughout the paper: straight-ahead components of a
$4$--regular projection, and the loop/monogon patterns corresponding to immediate Reidemeister~I reductions.

\medskip
\noindent
\emph{Straight-ahead components.}
For a $4$--regular map, each $\sigma$-cycle has length~$4$, and $\sigma^2$ pairs opposite darts at each vertex.
Define
\[
\rho:=\sigma^2\alpha.
\]

\begin{lemma}\label{lem:straight_ahead_components}
Let $P$ be a $4$--regular map encoded by $(\alpha,\sigma)$ and put
\[
\rho=\sigma^2\alpha.
\]
Then the number of straight-ahead components of the corresponding projection is
\[
c(P)=\frac12\,\#\mathrm{cyc}(\rho).
\]
\end{lemma}

\begin{proof}
A straight-ahead traversal alternates between crossing an edge, encoded by $\alpha$, and passing through a vertex to the opposite
dart, encoded by $\sigma^2$. Thus the oriented straight-ahead traversals are exactly the orbits of $\rho=\sigma^2\alpha$.
Each unoriented geometric component contributes two such orbits, one for each direction of traversal, giving the factor $1/2$.

In the closely related setting of immersed circles, Coquereaux--Zuber impose one-componentness by the condition that the corresponding
permutation $\sigma^2\tau$ consists of two cycles~\cite[Theorem~1]{CoquereauxZuber2016}.  In the present notation this is the special
case $\#\mathrm{cyc}(\rho)=2$, i.e. $c(P)=1$.
\end{proof}

\medskip
\noindent
\emph{Loop edges and monogons.}
An edge is a \emph{loop} (it joins a vertex to itself) precisely when its two darts lie in the same $\sigma$-cycle.
In the permutation model this is the condition that for some $h\in H$ the darts $h$ and $\alpha(h)$ belong to the same
vertex cycle of $\sigma$.

A monogon face is the projection-level shadow of an immediate Reidemeister~I reduction.
In the permutation model, monogons are detected by fixed points of the face permutation.

\begin{remark}[Monogons in the permutation model]\label{lem:monogon_fixed_point}
A face of degree~$1$ occurs if and only if the face permutation $\varphi=\sigma\alpha$ has a fixed point: by
\eqref{eq:classical_map_counts}, faces are the cycles of $\varphi$, and the degree of a face is the length of the corresponding cycle.
\end{remark}

The class “without planar loops” studied in \cite{KraskoOmelchenko2021NoPlanarLoops} is tailored to this
Reidemeister~I phenomenon: it excludes exactly those loop edges that bound a monogon face.
In the genus-one prime-projection setting used here, however, loop edges are also excluded for a graph-theoretic
reason: in a connected $4$--regular projection with at least two vertices, a loop edge forces a $2$--edge--cut.
Indeed, the loop uses two darts at its vertex, and the two remaining incident edges form a cut pair separating
that vertex together with the loop from the rest of the graph. Thus loop edges are ruled out by the same
projection-level obstruction that underlies Definition~\ref{def:composite_witness}.
The reference implementation nevertheless keeps looplessness as an explicit switch, so that larger, non-prime
projection pools can be reproduced if desired (Remark~\ref{rem:proj_optional_switches}).

\section{Prime projections on the torus: unsensed maps and canonical enumeration}
\label{sec:projections}

This section isolates the \emph{projection stage} of the genus-one tabulation pipeline.
The permutation language used here is the classical map language recalled in Section~\ref{sec:maps}; the new point is the way it is
specialized to a fixed-thickened-torus knot/link census.
Throughout this section, a \emph{projection} on the torus $T^2=\Sigma_1$ means a connected cellular embedding of a
$4$--valent multigraph, without over/under data; equivalently, a connected $4$--regular map on $T^2$.
The basic task is to enumerate such projections up to unsensed homeomorphism of the pair $(T^2,G)$, without omission
or duplication.
Unlike the immersed-circle setting of \cite{CoquereauxZuber2016}, we do not impose one-componentness of the straight-ahead curve.
Instead, the number of straight-ahead components is computed after the torus map has been generated, so that knots and links enter
the same projection pool.

As discussed in the Introduction, embedding-driven pipelines first enumerate abstract graphs and then classify their
embeddings into $T^2$, with the component structure becoming visible only after an embedding has been fixed.
In the map framework, unsensed equivalence becomes a finite orbit relation on permutation pairs
(Definition~\ref{prop:maps_conjugacy_equivalence}), canonical representatives give deterministic deduplication, and
straight-ahead components are read off from the same encoding
(Lemma~\ref{lem:straight_ahead_components}).
Thus knots and links can be stratified after a single projection pool has been generated.
We now formalize this projection stage.

\subsection{Projection model and standing constraints}
\label{subsec:proj_model_constraints}

Fix $N\ge 1$ and let $H=\{1,\dots,4N\}$ be the dart set.
We use the permutation encoding $(\alpha,\sigma)$ recalled in Section~\ref{sec:maps}
(with the same composition convention).  As usual,
$\sigma$ records the cyclic order of darts around vertices, $\alpha$ pairs darts into edges,
and $\varphi=\sigma\alpha$ is the face permutation.

For generation it is convenient to work with a fixed standard representative of the vertex rotation:
\begin{equation}\label{eq:standard_sigma_proj}
\sigma_0=(1\,2\,3\,4)(5\,6\,7\,8)\cdots(4N-3\ \ 4N-2\ \ 4N-1\ \ 4N).
\end{equation}
This does \emph{not} restrict generality: it only chooses a convenient labelling convention.

\begin{remark}[Standardization of the vertex rotation]\label{lem:standardize_sigma}
By relabelling darts, every labelled $4$--regular map with $N$ vertices can be encoded with vertex rotation
\[
\sigma_0=(1\,2\,3\,4)(5\,6\,7\,8)\cdots(4N-3\ \ 4N-2\ \ 4N-1\ \ 4N).
\]
Indeed, one orders the $N$ vertex cycles and chooses a starting dart in each cycle.  This standardization is only a choice of labels,
not a restriction on the unlabelled map.  The same normalization is used in closely related permutation enumerations of immersed
curves~\cite[Theorem~1]{CoquereauxZuber2016}.
\end{remark}

On the torus, Euler's formula for a $4$--regular map with $N$ vertices ($V=N$, $E=2N$) reads
\[
V-E+F=0 \qquad\Longleftrightarrow\qquad F=\#\mathrm{cyc}(\sigma\alpha)=N.
\]

\begin{definition}[Candidate torus projection]\label{def:candidate_torus_projection}
A \emph{labelled candidate projection on $T^2$ with $N$ vertices} is a pair $(\alpha,\sigma_0)$ such that:
\begin{enumerate}
\item $\alpha$ is a fixed-point-free involution on $H$;
\item \textbf{connectedness:} the action of $\langle\alpha,\sigma_0\rangle$ on $H$ is transitive;
\item \textbf{torus condition:} $\#\mathrm{cyc}(\sigma_0\alpha)=N$;
\item \textbf{no monogons:} $\varphi=\sigma_0\alpha$ has no fixed points;
\item \textbf{looplessness:} the projection has no loop edges, i.e.\ for every $h\in H$ the darts $h$ and $\alpha(h)$
      lie in different $\sigma_0$--cycles.
\end{enumerate}
An \emph{unlabelled candidate projection} is the corresponding unlabelled map, i.e.\ a homeomorphism class of pairs $(T^2,G)$,
as formalized in Definition~\ref{prop:maps_conjugacy_equivalence}.
\end{definition}

\begin{remark}[Monogons and looplessness]
By Remark~\ref{lem:monogon_fixed_point}, fixed points of $\varphi=\sigma\alpha$ are exactly monogon faces, which support
immediate Reidemeister~I reductions after crossing assignments.
Under the looplessness condition in Definition~\ref{def:candidate_torus_projection}, the no-monogon condition is in fact
redundant: if $\varphi(h)=h$, then $\alpha(\sigma_0(h))=h$, hence $\alpha(h)=\sigma_0(h)$, so the two darts of the
corresponding edge lie in the same $\sigma_0$--cycle.
We nevertheless state the no-monogon condition explicitly, since the framework and the reference implementation also allow
runs with loop edges, in which it remains a separate restriction.
\end{remark}

\begin{remark}[Optional projection-level switches]\label{rem:proj_optional_switches}
The map-theoretic framework itself does not depend on several minor convention choices, which we therefore state explicitly.

\smallskip
\noindent
\emph{(i) Bigons.}
A bigon face corresponds to a $2$--cycle of $\varphi=\sigma_0\alpha$.
We do \emph{not} exclude bigons at the projection level, since reduced or prime torus diagrams may still be supported by
projections with bigon faces.
When matching ``reduced diagram'' conventions, bigons are handled later by a diagram-level rule on crossing assignments;
see Definition~\ref{def:bigon_rule}.

\smallskip
\noindent
\emph{(ii) Allowing loop edges.}
Condition~(5) in Definition~\ref{def:candidate_torus_projection} excludes loop edges, whether or not they bound monogon faces.
This looplessness restriction is stronger than the planar-loop exclusion tailored to Reidemeister~I, but it is not merely
a convention: for connected $4$--regular projections with at least two vertices, any loop edge forces a $2$--edge--cut.
Such a projection would be discarded by Definition~\ref{def:composite_witness} in any case.
The reference implementation keeps looplessness as an explicit switch, so that the effect of allowing loop edges can be reproduced if desired.
\end{remark}

\subsection{Unsensed equivalence and canonical enumeration}
\label{subsec:proj_canonicalization}

Throughout the enumeration we work in the \emph{unsensed} regime, i.e.\ up to all
homeomorphisms of the pair $(T^2,G)$, including orientation reversal.
By Definition~\ref{prop:maps_conjugacy_equivalence}, two labelled encodings
$(\alpha_1,\sigma_1)$ and $(\alpha_2,\sigma_2)$ represent the same unlabelled
projection in this sense if and only if there exists $\pi\in S_H$ such that
\[
(\alpha_2,\sigma_2)=(\pi\alpha_1\pi^{-1},\,\pi\sigma_1\pi^{-1})
\quad\text{or}\quad
(\alpha_2,\sigma_2)=(\pi\alpha_1\pi^{-1},\,\pi\sigma_1^{-1}\pi^{-1}).
\]
This is a finite orbit relation on permutation pairs.
A convenient way to quotient it is to choose, in each orbit, a \emph{canonical representative}:
once every generated encoding is replaced by its canonical form, ordinary set deduplication
yields exactly one representative of each unlabelled projection class.

Fix a total order $\prec$ on labelled pairs $(\alpha,\sigma)$ (e.g.\ lex order on the list of images of $1,\dots,|H|$
under $\alpha$ and $\sigma$).
For a labelled encoding $(\alpha,\sigma)$ let
\[
\mathcal{O}_{\mathrm{sensed}}(\alpha,\sigma)=\{(\pi\alpha\pi^{-1},\pi\sigma\pi^{-1}) : \pi\in S_H\},
\]
and let $\mathcal{O}_{\mathrm{unsensed}}(\alpha,\sigma)$ be the union of
$\mathcal{O}_{\mathrm{sensed}}(\alpha,\sigma)$ and $\mathcal{O}_{\mathrm{sensed}}(\alpha,\sigma^{-1})$.

\begin{definition}[Canonical representative]\label{def:proj_canonical_rep}
The \emph{sensed} canonical representative of $(\alpha,\sigma)$ is the $\prec$--minimal element of
$\mathcal{O}_{\mathrm{sensed}}(\alpha,\sigma)$.
The \emph{unsensed} canonical representative is the $\prec$--minimal element of
$\mathcal{O}_{\mathrm{unsensed}}(\alpha,\sigma)$.
\end{definition}

Each finite orbit has a unique $\prec$--minimal element.  Therefore two labelled encodings define the same unlabelled map in the
sensed or unsensed sense if and only if their corresponding canonical representatives coincide.
Consequently, once each produced encoding is replaced by its canonical representative, ordinary set deduplication leaves
\emph{exactly one} representative of each unlabelled projection class.\label{lem:proj_canonical_complete}
In practice, the same unsensed canonical representative is computed via a rooted normalization procedure, avoiding
enumeration of the full group $S_H$; see Appendix~\ref{app:proj_algorithms}.

We now formalize the enumeration step itself.
By Remark~\ref{lem:standardize_sigma}, it suffices to enumerate labelled encodings of the form $(\alpha,\sigma_0)$.

\medskip
\noindent\textbf{Conceptual enumeration algorithm.}
Enumerate all fixed-point-free involutions $\alpha$ on $H$ (equivalently, all perfect matchings of $H$).
For each $\alpha$:
\begin{enumerate}
\item test the defining conditions of Definition~\ref{def:candidate_torus_projection};
\item if the conditions hold, compute the unsensed canonical representative of $(\alpha,\sigma_0)$ and store it.
\end{enumerate}
Finally, deduplicate the stored canonical pairs.

This description is intentionally implementation-free: it isolates the mathematical content.
Pruning strategies that reduce the number of matchings explored (canonical growth rules, early rejection, etc.)
are efficiency devices and do not change the abstract correctness statement below.

\begin{theorem}[Correctness of projection enumeration]\label{thm:proj_enumeration_correctness}
Fix $N\ge 1$.  After deduplication, the procedure above outputs a set $\mathcal{C}_N$ of unsensed canonical representatives
with the following properties:
\begin{enumerate}
\item \textbf{Soundness:} every element of $\mathcal{C}_N$ represents an unlabelled candidate torus projection
      (Definition~\ref{def:candidate_torus_projection});
\item \textbf{Completeness:} every unlabelled candidate torus projection with $N$ vertices occurs (up to unsensed
      equivalence) among the outputs;
\item \textbf{Non-duplication:} no two distinct elements of $\mathcal{C}_N$ represent the same unsensed projection class.
\end{enumerate}
\end{theorem}

\begin{proof}
Soundness is immediate from the tests in Definition~\ref{def:candidate_torus_projection}: every stored pair is first checked to be a
connected, loopless, monogon-free torus projection, and canonicalization does not change its unlabelled class.

For completeness, let $P$ be any unlabelled candidate torus projection with $N$ vertices.  Choose an orientation and dart labels.
By the standardization in Remark~\ref{lem:standardize_sigma}, the resulting labelled encoding may be relabelled to the form
$(\alpha',\sigma_0)$.  The perfect matching $\alpha'$ is among those enumerated, and it passes the defining tests because it encodes $P$.
Thus the unsensed canonical representative of the class of $P$ is stored.

Finally, two stored outputs representing the same unsensed class have the same unsensed canonical representative by the definition of
canonical representative.  Removing exact duplicates therefore leaves exactly one representative of each unsensed class.
\end{proof}

\begin{remark}[Efficiency and implementation]\label{rem:efficiency}
Theorem~\ref{thm:proj_enumeration_correctness} is independent of how perfect matchings
are generated.
In practice we construct $\alpha$ recursively, applying early rejection
(connectedness, torus face-count, monogon/loop constraints) as soon as they become forced,
and use a canonical construction path (``activated vertices'') to suppress relabelling
redundancy compatible with~$\sigma_0$.
Unsensed canonical representatives are computed by rooted normalization of the dart graph
over all root darts and both orientations $\sigma$~vs.~$\sigma^{-1}$, avoiding enumeration
of the full group~$S_H$.
These are purely efficiency devices; a detailed description is given in
Appendix~\ref{app:proj_algorithms}.
\end{remark}

\subsection{Components and projection-level primeness witnesses}
\label{subsec:proj_components}

For a labelled encoding $(\alpha,\sigma)$ of a $4$--regular projection, set
\[
\rho(\alpha,\sigma):=\sigma^2\alpha .
\]
By Lemma~\ref{lem:straight_ahead_components},
\begin{equation}\label{eq:component_count_torus_proj}
c(P)=\tfrac12\,\#\mathrm{cyc}(\rho(\alpha,\sigma)).
\end{equation}
In particular, for the standardized encodings $(\alpha,\sigma_0)$ used during generation this becomes
\[
c(P)=\tfrac12\,\#\mathrm{cyc}(\sigma_0^2\alpha).
\]
Thus, once a representative encoding of each class has been fixed, the set $\mathcal{C}_N$ can be stratified into knot
projections ($c(P)=1$) and link projections ($c(P)\ge 2$).

To match the notion of ``prime projections'' used in genus-one tabulations, we apply two projection-level filters
that depend only on the embedded graph (hence only on $(\alpha,\sigma_0)$) and are therefore stable under the
equivalences of Definition~\ref{prop:maps_conjugacy_equivalence}.
These filters are stated explicitly as computable witnesses and are used to align with published genus-one projection
tables; they are not intended as an exhaustive topological primeness criterion.

Let $(\alpha,\sigma)$ be any labelled encoding of $P$.
The underlying abstract multigraph $G(P)$ has the $\sigma$-cycles as vertices and the $\alpha$-pairs as edges; for the
standardized encodings used during generation, these are the $\sigma_0$-cycles.

\begin{definition}[Compositeness witness]\label{def:composite_witness}
We call $P$ \emph{$2$--edge--cut composite} if the multigraph $G(P)$ has a $2$--edge--cut, i.e.\ a set of two edges whose
removal disconnects $G(P)$.
\end{definition}

\begin{remark}[Disc-separating witness]\label{lem:connected_sum_implies_2edgecut}
If a simple closed curve $C\subset T^2$ bounds an embedded disc and meets the projection transversely in exactly two interior
points of edges, with vertices on both sides of $C$, then the two crossed edges are the only edges joining the inside part to the
outside part.  Hence they form a $2$--edge--cut in $G(P)$.
\end{remark}

For links we also use a fast splitness witness based on mixed vertices (mixed crossings in the diagrammatic language).
Compute straight-ahead components using \eqref{eq:component_count_torus_proj}, and call a vertex (that is, a $\sigma$-cycle)
\emph{mixed} if its four darts do not all lie in the same straight-ahead component.

\begin{definition}[Splitness witness]\label{def:split_witness}
For a link projection $P$ with $c(P)\ge 2$, and for a straight-ahead component $K$ of $P$, let $m(K)$ be the number of
mixed vertices incident to at least one dart of $K$.
We say that $P$ is \emph{split-witnessed} if there exists a component $K$ with $m(K)\le 1$.
\end{definition}

\begin{definition}[Prime projections used in this paper]\label{def:prime_projection_computational}
A candidate torus projection $P$ is called \emph{prime} (for the purposes of this paper) if:
\begin{itemize}
\item if $c(P)=1$ then $P$ is not $2$--edge--cut composite;
\item if $c(P)\ge 2$ then $P$ is not $2$--edge--cut composite and is not split-witnessed.
\end{itemize}
\end{definition}

\begin{remark}[Scope of the definition]
Definition~\ref{def:prime_projection_computational} is intentionally phrased in terms of explicit witnesses that can be
checked directly from the projection encoding and are stable under unsensed equivalence.
Remark~\ref{lem:connected_sum_implies_2edgecut} shows that the $2$--edge--cut condition is a necessary obstruction to
connected-sum decompositions.

For the splitness witness, note first that $m(K)=0$ cannot occur in a connected link projection:
if every vertex incident to $K$ were unmixed, then all edges incident to those vertices would belong to $K$, so $K$
would form a connected component of $G(P)$.
Thus, in the connected setting, the witness concerns the case $m(K)=1$, where one straight-ahead component meets the
remainder of the projection at a single mixed vertex.
We use this as an explicit projection-level filter to match the convention adopted in the genus-one tables of
\cite{AkimovaMatveevTarkaev2020}, and we state it separately so that no hidden assumptions enter later comparisons.
\end{remark}

\subsection{Genus-one benchmark data}
\label{subsec:g1_projection_benchmark}

As a practical benchmark for Theorem~\ref{thm:proj_enumeration_correctness}, we record the size of the sets of unsensed
projection classes and the effect of the primeness filters for small $N$.
The numbers in Tables~\ref{tab:g1_projection_counts} and~\ref{tab:g1_links_by_components} are obtained by running the
enumeration procedure and counting unsensed canonical representatives.
For each $N\le 10$, the reference code~\cite{OmelchenkoTorusMaps} provides machine-readable lists of
unsensed canonical representatives for both the candidate and prime projection pools, so that all
counts and stratifications reported here can be regenerated directly from the released data.

\begin{table}[t]
\centering
\begin{tabular}{c|ccc|cc}
\toprule
$N$ & candidates & removed by $2$--edge--cut & prime total & knot projections & link projections \\
\midrule
$3$  & $6$       & $0$       & $6$       & $2$       & $4$       \\
$4$  & $28$      & $5$       & $23$      & $10$      & $13$      \\
$5$  & $109$     & $28$      & $81$      & $34$      & $47$      \\
$6$  & $595$     & $216$     & $379$     & $170$     & $209$     \\
$7$  & $3216$    & $1421$    & $1795$    & $777$     & $1018$    \\
$8$  & $19956$   & $10141$   & $9815$    & $4308$    & $5507$    \\
$9$  & $127085$  & $71492$   & $55593$   & $23783$   & $31810$   \\
$10$ & $845961$  & $514728$  & $331233$  & $139780$  & $191453$  \\
\bottomrule
\end{tabular}
\caption{Unsensed projection counts on the torus $T^2$ for $N=3,\ldots,10$.
The column ``candidates'' counts unlabelled candidate projections up to unsensed equivalence in the sense of
Definition~\ref{def:candidate_torus_projection}.
The column ``removed by $2$--edge--cut'' records how many candidates are discarded by the compositeness witness
(Definition~\ref{def:composite_witness}); in the range $N\le 10$, the splitness witness
(Definition~\ref{def:split_witness}) removes no further candidates.
The remaining prime projections are stratified by the number of straight-ahead components into knot projections
($c(P)=1$) and link projections ($c(P)\ge 2$).}
\label{tab:g1_projection_counts}
\end{table}

In the range $N\le 10$ under the standing constraints of Definition~\ref{def:candidate_torus_projection}, the splitness
witness does not remove additional projections beyond the $2$--edge--cut filter.
Accordingly, Table~\ref{tab:g1_projection_counts} records only the effect of the $2$--edge--cut witness.

\begin{table}[t]
\centering
\begin{tabular}{r|cccccc}
\toprule
$N$ & $c(P)=2$ & $c(P)=3$ & $c(P)=4$ & $c(P)=5$ & $c(P)=6$ & $c(P)=7$ \\
\midrule
$3$  & $3$      & $1$     & $0$    & $0$   & $0$  & $0$ \\
$4$  & $9$      & $3$     & $1$    & $0$   & $0$  & $0$ \\
$5$  & $37$     & $9$     & $1$    & $0$   & $0$  & $0$ \\
$6$  & $150$    & $51$    & $7$    & $1$   & $0$  & $0$ \\
$7$  & $775$    & $212$   & $30$   & $1$   & $0$  & $0$ \\
$8$  & $4030$   & $1293$  & $169$  & $14$  & $1$  & $0$ \\
$9$  & $23484$  & $7276$  & $994$  & $54$  & $2$  & $0$ \\
$10$ & $138377$ & $46112$ & $6502$ & $444$ & $17$ & $1$ \\
\bottomrule
\end{tabular}
\caption{Distribution of unsensed prime \emph{link} projections on the torus $T^2$ by the number of straight-ahead
components $c(P)\ge 2$.}
\label{tab:g1_links_by_components}
\end{table}

\begin{remark}[Higher genus]
The proof of Theorem~\ref{thm:proj_enumeration_correctness} uses only the permutation model
and canonical representatives; for $\Sigma_g$ the genus condition becomes
$\#\mathrm{cyc}(\sigma\alpha)=N+2-2g$, and the argument applies verbatim.
The additional complications that arise at the diagram level for $g\ge 2$ are discussed in
Section~\ref{sec:outlook}.
\end{remark}

\section{From projections to diagrams on the torus: crossing data, local reductions, and state-sum invariants}
\label{sec:diagrams}

Starting from the pool of unsensed prime \emph{projections} on the torus constructed in
Section~\ref{sec:projections}, we now pass to \emph{diagrams}, i.e.\ projections endowed with
over/under information at crossings.
For a fixed projection, this additional datum is purely local: at each vertex one chooses one of the two local strands
as the overpassing strand.
After fixing a standardized representative of the projection, this choice is encoded by one bit per vertex.

The purpose of this section is to keep the diagram stage intrinsic.
We record the diagram-level conventions used in genus-one tabulations, prove a local
Reidemeister~II criterion supported by a bigon face, and formulate the generalized
Kauffman-type bracket entirely within the permutation model.
In particular, the state circles of a smoothing arise from explicit permutations, and the
contractible/essential dichotomy on $T^2$ is detected by a finite $\mathbb{F}_2$--homology test.
No geometric drawing in a fundamental polygon is required at any stage.

\subsection{Crossing assignments and local diagram-level rules}
\label{subsec:diag_assignments_and_rules}

Let $P$ be a fixed projection on $T^2$, encoded by a pair $(\alpha,\sigma)$ as in
Section~\ref{sec:maps}. Its vertices are the $\sigma$-cycles, and at each vertex the two local strands
are the orbits of the opposite-dart involution $\sigma^2$.

\begin{definition}[Diagram supported by a projection]\label{def:diagram_on_projection}
A \emph{crossing assignment} on $P$ chooses, for each vertex $v$, one of the two $\sigma^2$-orbits at $v$ as the
overpassing local strand.
The resulting \emph{diagram} is denoted by $D=(P,b)$.
\end{definition}

For computation, we work with a standardized representative $(\alpha,\sigma_0)$ as in
Section~\ref{sec:projections}. If the $k$th vertex is the block
\[
(4k-3\,\,4k-2\,\,4k-1\,\,4k),
\]
we encode the crossing assignment by a bit $b(v_k)\in\{0,1\}$:
$b(v_k)=0$ means that the strand $\{4k-3,4k-1\}$ is overpassing, and
$b(v_k)=1$ means that the strand $\{4k-2,4k\}$ is overpassing.

More generally, whenever a local cyclic order
\[
(h_0\,h_1\,h_2\,h_3)
\]
is chosen at a vertex $v$, we encode the crossing choice at $v$ by a bit relative to that order:
$b(v)=0$ if the local strand $\{h_0,h_2\}$ is overpassing, and $b(v)=1$ if the local strand $\{h_1,h_3\}$ is overpassing.
For the standardized representative, this agrees with the blockwise convention above.

Thus, a projection with $N$ crossings supports $2^N$ diagrams at the purely combinatorial level.
In what follows we often treat $P$ as fixed and vary $b$.

\begin{remark}[Global crossing switch as an optional convention]\label{rem:global_switch}
In the standardized bit encoding, the involution $b\mapsto 1-b$ reverses the over/under choice at every crossing
simultaneously.
Some genus-one tabulations quotient by this symmetry when organizing lists of diagrams on a fixed projection.
Since this is a matter of convention (mirror-related diagrams may or may not be identified), we treat it as an explicit
switch: when enabled, we enumerate crossing assignments modulo $b\sim 1-b$ by fixing one bit (e.g.\ $b(v_1)=0$).
\end{remark}

\paragraph{Bigon reduction criterion.}

We next explain how an immediate Reidemeister~II reduction appears in the permutation model.
This criterion is intrinsic: it expresses a genuine local Reidemeister~II reduction on the surface and does not depend on
table-matching conventions.

\begin{definition}[Bigon faces]\label{def:bigon_face}
Let $P=(\alpha,\sigma)$ be a projection and let $\varphi=\sigma\alpha$ be its face permutation.
A \emph{bigon face} is a face of degree~$2$, equivalently a $2$--cycle $(i\,j)$ of $\varphi$.
Since $P$ is cellular, every face (hence every bigon face) is a topological disc.
\end{definition}

Write $(i\,j)$ for a $2$--cycle of $\varphi$.  By definition of $\varphi=\sigma\alpha$ and our composition convention,
\[
\varphi(i)=j \iff \alpha(\sigma(i))=j,\qquad
\varphi(j)=i \iff \alpha(\sigma(j))=i.
\]
Thus the bigon boundary uses the two edges
\[
e_1=\{\sigma(i),j\},\qquad e_2=\{\sigma(j),i\},
\]
connecting the two incident vertices.

\begin{lemma}[Bigon $\Omega_2$ criterion in crossing bits]\label{lem:bigon_omega2}
Let $P$ be a projection, let $(i\,j)$ be a $2$--cycle of $\varphi=\sigma\alpha$ incident to two distinct
vertices $u\neq v$, and let $D=(P,b)$ be a diagram supported by $P$.
Choose starting darts and write the $\sigma$-cycles of $u$ and $v$ as
\[
(h_0\,h_1\,h_2\,h_3),\qquad (k_0\,k_1\,k_2\,k_3),
\]
with $\sigma(h_r)=h_{r+1}$ and $\sigma(k_s)=k_{s+1}$.
Encode the crossing assignment at $u$ and $v$ by bits relative to these local orders:
$b(u)=0$ if $\{h_0,h_2\}$ is overpassing and $b(u)=1$ if $\{h_1,h_3\}$ is overpassing, and similarly for $b(v)$.
If $i=h_r$ and $j=k_s$, define the local parities
\[
p_u(i):=r \pmod 2,\qquad p_v(j):=s \pmod 2.
\]
Then $D$ admits an immediate Reidemeister~II reduction supported in the bigon face corresponding to $(i\,j)$
if and only if
\begin{equation}\label{eq:bigon_rule_bits}
b(u)+b(v)+p_u(i)+p_v(j)\equiv 1 \pmod 2.
\end{equation}
Moreover, the congruence \eqref{eq:bigon_rule_bits} is independent of the chosen starting darts in the two $\sigma$-cycles.
\end{lemma}

\begin{proof}
Since the bigon is a disc face, a neighbourhood of its boundary is homeomorphic to the standard planar
Reidemeister~II neighbourhood.

At $u$, the bigon uses the adjacent darts $i$ and $\sigma(i)$, hence it meets the two local strands at $u$ once each.
In the chosen local encoding, the overpassing bigon edge at $u$ is the edge incident to $i$ if and only if
$b(u)=p_u(i)$, and otherwise it is the edge incident to $\sigma(i)$.
Similarly, at $v$ the overpassing bigon edge is incident to $j$ if and only if $b(v)=p_v(j)$.

Within the bigon neighbourhood, the two geometric strands pair the incident darts as
\[
(i,\sigma(j))\qquad\text{and}\qquad (\sigma(i),j).
\]
Hence an immediate $\Omega_2$ cancellation occurs exactly when the same geometric strand is overpassing at both vertices.
This happens precisely in the two cases
\[
b(u)=p_u(i),\quad b(v)=1-p_v(j),
\]
or
\[
b(u)=1-p_u(i),\quad b(v)=p_v(j),
\]
which is equivalent to \eqref{eq:bigon_rule_bits}.

Finally, changing the chosen starting dart at $u$ by one step interchanges the two local bit values and also flips
$p_u(i)$ modulo~$2$; thus $b(u)+p_u(i)$ is unchanged modulo~$2$. The same holds at $v$.
Therefore the congruence \eqref{eq:bigon_rule_bits} is independent of the chosen starting darts.
\end{proof}

\begin{definition}[Bigon rule used for reduced-diagram conventions]\label{def:bigon_rule}
When we adopt the reduced-diagram convention used in the genus-one $N\le 4$ tables of \cite{AkimovaMatveev2012},
we keep only those crossing assignments $b$ for which no two-vertex bigon supports an immediate $\Omega_2$ reduction.
By Lemma~\ref{lem:bigon_omega2}, this is equivalent to requiring
\[
b(u)+b(v)+p_u(i)+p_v(j)\equiv 0 \pmod 2
\]
for every bigon face $(i\,j)$ incident to two distinct vertices $u\neq v$, where the local bits and parities are computed
with respect to any chosen starting darts in the two incident $\sigma$-cycles; the condition is independent of that choice.
\end{definition}

\begin{remark}[Classical special case]
When $p_u(i)=p_v(j)=0$ (for example, when both $i$ and $j$ are the first darts of their standardized blocks),
the criterion~\eqref{eq:bigon_rule_bits} reduces to $b(u)\neq b(v)$, which is the familiar naive form in local planar
labellings.
\end{remark}

\begin{remark}[Degenerate bigons incident twice to one vertex]
If a bigon face is incident twice to the same vertex ($u=v$), then the local picture is not the standard two-vertex
Reidemeister~II configuration.  In our computations this case is treated separately (and does not contribute to the
rule above).
\end{remark}

\subsection{The generalized bracket in permutation form}
\label{subsec:diag_bracket}

We now turn to invariants.
For diagrams on the torus we use the generalized Kauffman-type bracket introduced in
\cite{AkimovaMatveev2012,AkimovaMatveevTarkaev2018}.
We recall the state-sum definition and then prove that every ingredient of the state sum can be obtained directly from
the permutation encoding $(\alpha,\sigma)$ and a crossing assignment $b$.

\medskip
\noindent\textbf{States, smoothings, and the state permutation.}

Fix a diagram $D=(P,b)$ where $P=(\alpha,\sigma)$.
To describe smoothings combinatorially, we first build the involutions that encode the smoothing connections
\emph{inside vertices}, and then define states as choices among these smoothings at each crossing.

For each vertex $v$, let $O_v$ be the overpassing local strand selected by the crossing assignment $b$, and
recall the opposite-dart involution $\sigma^2$ (which swaps the two darts of each local strand at $v$).
Choose any dart $a_v\in O_v$ and define two pairings of the four darts at $v$ by
\begin{align*}
\tau_v^{A} &:= (a_v\,\sigma(a_v))(\sigma^2(a_v)\,\sigma^3(a_v)),\\
\tau_v^{B} &:= (\sigma(a_v)\,\sigma^2(a_v))(\sigma^3(a_v)\,a_v).
\end{align*}
These pairings are independent of the choice of $a_v\in O_v$: replacing $a_v$ by $\sigma^2(a_v)$ merely exchanges the
two transpositions in each product.
They are the $A$-- and $B$--smoothings at the crossing (in the sense of the Kauffman bracket), labelled relative to the
overpassing strand $O_v$.
If the overpassing strand is switched to the complementary local strand, then $\tau_v^{A}$ and $\tau_v^{B}$ are interchanged.

A \emph{state} is a function
\[
s:\mathrm{Vert}(P)\longrightarrow \{A,B\},
\]
and we write $n_A(s)$ and $n_B(s)$ for the numbers of vertices where $s$ chooses the $A$-- and $B$--smoothing,
respectively.

For a state $s$, define the global smoothing involution $\tau_{b,s}$ on $H$ by taking the product of
$\tau_v^{s(v)}$ over all vertices $v$.

Finally, define the \emph{state permutation} by
\begin{equation}\label{eq:pi_bs_def_section}
\pi_{b,s}:=\alpha\,\tau_{b,s}.
\end{equation}
Recall that in our convention $(pq)(h)=q(p(h))$, so $\pi_{b,s}(h)=\tau_{b,s}(\alpha(h))$:
one traverses an edge (apply $\alpha$) and then follows the smoothing inside the next vertex (apply $\tau_{b,s}$).

\begin{lemma}[State circles from $(\alpha,\tau_{b,s})$]\label{lem:state_circles_cycles_section}
Fix a diagram $D=(P,b)$ and a state $s$, and let $\tau_{b,s}$ be the smoothing involution.
Consider the graph on the dart set $H$ whose edges are the $\alpha$--pairs and the $\tau_{b,s}$--pairs.
Each connected component of this graph is a cycle, and these components are in canonical bijection with the
(state) circles obtained by smoothing $D$ according to $s$.

Moreover, $\pi_{b,s}=\alpha\,\tau_{b,s}$ advances by two steps along each such dart-cycle.
Consequently, each geometric state circle gives rise to exactly two cycles of $\pi_{b,s}$ (the two orientations), and
\[
\#\{\text{state circles of the smoothing }(D,s)\}=\tfrac12\,\#\mathrm{cyc}(\pi_{b,s}).
\]
\end{lemma}

\begin{proof}
The involutions $\alpha$ and $\tau_{b,s}$ pair darts in the two evident ways: $\alpha$ pairs the two darts of each edge,
and $\tau_{b,s}$ pairs the darts connected by the chosen smoothing inside each vertex.
Hence the graph with $\alpha$-- and $\tau_{b,s}$--edges is the combinatorial dart model of the smoothed diagram.
Every dart is incident to exactly one $\alpha$--edge and exactly one $\tau_{b,s}$--edge, so each connected component is
a $2$--regular graph, hence a cycle.
Because the two kinds of edges alternate, every such cycle has even length, and these cycles are precisely the state circles
of the smoothing.

Finally, one step along an oriented smoothed circle traverses an edge (apply $\alpha$) and then follows the smoothing
inside the next vertex (apply $\tau_{b,s}$), which is exactly $\pi_{b,s}=\alpha\tau_{b,s}$ in our right-action convention.
If a dart-cycle has length $2m$, then $\pi_{b,s}$ acts on it by moving two steps at a time, hence has exactly
$\gcd(2,2m)=2$ orbits on that cycle.
These two $\pi_{b,s}$--cycles correspond to the two orientations of the same geometric state circle.
\end{proof}

\medskip
\noindent\textbf{Contractible vs.\ essential state circles.}

The bracket on the torus distinguishes contractible and essential state circles.
We now show that this distinction can be detected without any geometric drawing.

Let $E$ be the set of (undirected) edges of the projection, i.e.\ the set of $\alpha$-orbits in $H$.
Let $C_1$ be the $\mathbb{F}_2$--vector space freely generated by $E$.
The map $P$ gives a cellular decomposition of $T^2$, hence also a cellular chain complex
\[
C_2 \xrightarrow{\partial_2} C_1 \xrightarrow{\partial_1} C_0
\qquad (\text{coefficients in }\mathbb{F}_2).
\]
The image $\mathrm{im}(\partial_2)\subset C_1$ is the $\mathbb{F}_2$--subspace of $1$--boundaries.

Every state circle $C$ meets edges of the projection in a finite collection of transverse arcs.
Define $v(C)\in C_1$ to be the edge-incidence vector modulo~$2$: the coefficient of an edge $e\in E$ is the parity of the
number of times $C$ traverses $e$.

Similarly, each face $f$ of the projection determines its boundary vector $\partial_2(f)\in C_1$, whose coefficients are the
parities of edge traversals along the boundary walk of $f$.

\begin{lemma}[Face span equals the boundary subspace]\label{lem:face_span_boundaries_section}
The subspace $\mathrm{im}(\partial_2)\subset C_1$ is spanned by the face boundary vectors $\partial_2(f)$, where $f$ ranges
over faces of the projection.
\end{lemma}

\begin{proof}
This is the definition of the cellular boundary map $\partial_2$: it sends each $2$--cell (face) to its attaching
$1$--cycle.  Hence $\mathrm{im}(\partial_2)$ is generated by these images.
\end{proof}

\begin{lemma}[Contractible state circles on $T^2$ via $\mathbb{F}_2$]\label{lem:cut_via_face_span_section}
Let $C$ be a state circle of a torus diagram, viewed as an embedded simple closed curve on $T^2$, and let
$v(C)\in C_1$ be its edge-incidence vector mod~$2$.
Then the class $[C]$ is zero in $H_1(T^2;\mathbb{F}_2)$ if and only if $v(C)\in \mathrm{im}(\partial_2)$.
Moreover, this is equivalent to $C$ being contractible.
\end{lemma}

\begin{proof}
Choose pairwise disjoint closed discs $B_v$ around the vertices of the projection.
For a fixed state circle $C$, the intersection $C\cap B_v$ is a union of $0$, $1$, or $2$ smoothing arcs.
Each such smoothing arc is properly embedded in $B_v$ and joins two boundary points lying on incident half-edges.

Let $v(C)\in C_1$ be the cellular $1$--chain whose coefficient on an edge $e$ is the parity of the number of times
$C$ traverses $e$.
At each vertex, the smoothing pairs the incident half-edges, so an even number of half-edges of $v(C)$ are incident to
that vertex. Hence $\partial_1 v(C)=0$, and $v(C)$ is a cellular $1$--cycle.

Now fix one smoothing arc $\gamma\subset C\cap B_v$, and let $\eta_1,\eta_2$ be the two half-edge segments in $B_v$
joining the endpoints of $\gamma$ to the vertex $v$.
Then $\gamma\cup\eta_1\cup\eta_2$ bounds a small $2$--chain in $B_v$.
Over $\mathbb{F}_2$, replacing every smoothing arc of $C$ in this way changes $C$ by a boundary.
After carrying out this replacement in every vertex disc, one obtains precisely the cellular cycle $v(C)$.
Therefore
\[
[C]=[v(C)]\in H_1(T^2;\mathbb{F}_2).
\]
It follows that $[C]=0$ if and only if $v(C)$ is a cellular boundary, i.e.\ if and only if
$v(C)\in\mathrm{im}(\partial_2)$, proving the first statement.

For the second statement, every embedded essential simple closed curve on $T^2$ represents a primitive class in
$H_1(T^2;\mathbb{Z})\cong\mathbb{Z}^2$, say $(p,q)$ with $\gcd(p,q)=1$.
Such a class satisfies $(p,q)\not\equiv(0,0)\pmod 2$, so its reduction in $H_1(T^2;\mathbb{F}_2)$ is nonzero.
Therefore $[C]=0$ in $H_1(T^2;\mathbb{F}_2)$ if and only if $[C]=0$ in $H_1(T^2;\mathbb{Z})$, which for an embedded
simple closed curve is equivalent to $C$ being contractible.
\end{proof}

\medskip
\noindent\textbf{The state sum.}

For a state $s$, let $\gamma(s)$ and $\delta(s)$ denote the numbers of contractible and essential \emph{geometric} state circles.
By Lemma~\ref{lem:state_circles_cycles_section}, each geometric state circle contributes two cycles of $\pi_{b,s}$, corresponding
to the two orientations, and both orientations traverse the same set of edges.
Thus contractibility/essentiality is the same for the two corresponding $\pi_{b,s}$--cycles, and it is detected by the
$\mathbb{F}_2$ boundary test of Lemma~\ref{lem:cut_via_face_span_section}.

\begin{definition}[Generalized bracket on the torus]\label{def:generalized_bracket_section}
For a torus diagram $D=(P,b)$ with $N$ crossings, the generalized bracket is
\begin{equation}\label{eq:bracket_state_sum_section}
\langle D\rangle(a,x)=\sum_{s:\mathrm{Vert}(P)\to\{A,B\}}
a^{n_A(s)-n_B(s)}\bigl(-a^2-a^{-2}\bigr)^{\gamma(s)}x^{\delta(s)}.
\end{equation}
\end{definition}

\begin{remark}[Invariance]\label{rem:bracket_invariance}
The bracket \eqref{eq:bracket_state_sum_section} and its normalization discussed below coincide with the invariants
introduced and used in \cite{AkimovaMatveev2012,AkimovaMatveevTarkaev2018}.
In particular, $\langle D\rangle(a,x)$ is invariant under Reidemeister moves of types $\Omega_2$ and $\Omega_3$, and its
behaviour under $\Omega_1$ is the standard multiplicative factor familiar from the classical Kauffman bracket.
We refer to the cited sources for the full invariance proof and sign conventions.
\end{remark}

\paragraph{Normalization.}

For knots on the torus, one can normalize the bracket by writhe to obtain an invariant under $\Omega_1$ as well.

\begin{definition}[Writhe and $X_D(a,x)$ for knots]\label{def:Xpoly_section}
Let $D$ be a \emph{knot} diagram on $T^2$ (i.e.\ $c(P)=1$), and equip the knot with an orientation.
The writhe $w(D)$ is the sum of crossing signs in the usual sense, computed using the orientation of $T^2$.
For a knot, reversing the orientation does not change $w(D)$, so the following definition is orientation-independent:
\[
X_D(a,x)=(-a)^{-3w(D)}\langle D\rangle(a,x).
\]
\end{definition}

\begin{remark}[Links and orientation dependence]
For links ($c(P)\ge 2$) the writhe depends on a choice of component orientations, so $X_D(a,x)$ is not, by itself, an
invariant of an unoriented link.
In computations we therefore use the orientation-independent bracket $\langle D\rangle(a,x)$ as the primary classification
key for links, in accordance with \cite{AkimovaMatveevTarkaev2020}.
\end{remark}

\subsection{Skeletons and comparison conventions}
\label{subsec:diag_skeleton_conventions}

For knots we apply the following construction to $X_D(a,x)$, and for links to the bracket
$\langle D\rangle(a,x)$.
Write generically
\[
Q_D(a,x)=\sum_m P_m(a)\,x^m,\qquad P_m(a)\in\mathbb{Z}[a,a^{-1}],
\]
where
\[
Q_D=
\begin{cases}
X_D(a,x), & \text{for knots},\\[2mm]
\langle D\rangle(a,x), & \text{for links}.
\end{cases}
\]

For each $x$--degree $m$, write
\[
P_m(a)=\sum_k c_{m,k}a^k,
\qquad c_{m,k}\in\mathbb{Z},
\]
and let
\[
S_m(Q_D):=(c_{m,k_1},c_{m,k_2},\dots,c_{m,k_r}),
\]
where $k_1<k_2<\cdots<k_r$ are the exponents for which $c_{m,k}\neq 0$.
Thus $S_m(Q_D)$ records the nonzero coefficients of $P_m(a)$ in increasing exponent order and forgets the
corresponding $a$--exponents.

The \emph{skeleton} of $Q_D$ is the indexed family
\[
\operatorname{Sk}(Q_D):=(S_m(Q_D))_m.
\]
It retains the coefficient pattern in each $x$--degree while discarding the $a$--exponent data.

For comparison with published tables, we may further replace $\operatorname{Sk}(Q_D)$ by a canonical representative
under the transformations induced by $a\mapsto a^{-1}$ (which reverses each tuple $S_m$) and, when required,
by overall multiplication of $Q_D$ by $-1$ (which changes all entries of all tuples simultaneously).

\begin{remark}[Skeleton is not a complete invariant]
Even after the comparison-level canonicalizations just described, skeleton equality is only a coarse necessary condition
for equality of the full invariant under the same comparison conventions.
Distinct Laurent polynomials can share the same coefficient pattern while differing in their $a$--exponents.
Accordingly, the skeleton is used only for inexpensive grouping, and final separation uses the full bracket or polynomial.
\end{remark}

\paragraph{Over/under participation filter.}

A second convention that affects comparison with published link tables is the \emph{over/under participation
filter}, used in some genus-one link tabulations.
Unlike the bigon rule of \S\ref{subsec:diag_assignments_and_rules}, this is not a topological reduction statement;
it is an explicit convention that removes diagrams with a ``layered'' mixed-crossing pattern.

Let $P$ be a projection, and let $c(P)$ be the number of straight-ahead components (computed from $P$ alone, see
Lemma~\ref{lem:straight_ahead_components} in Section~\ref{sec:maps}).
A vertex $v$ is \emph{mixed} if its four darts do not all belong to the same straight-ahead component.

\begin{definition}[Over/under participation filter]\label{def:overunder_filter}
Let $D=(P,b)$ be a link diagram, i.e.\ $c(P)\ge 2$.
We say that $D$ satisfies the \emph{over/under participation condition} if for each straight-ahead component $K$ of $P$
there exists a mixed vertex at which the local strand belonging to $K$ is overpassing, and a (possibly different) mixed
vertex at which the local strand belonging to $K$ is underpassing.
\end{definition}

\begin{remark}[Convention status]\label{rem:overunder_convention}
The filter of Definition~\ref{def:overunder_filter} is used in \cite{AkimovaMatveevTarkaev2020} as a fast prefilter, but it
is not logically forced by primeness or non-splitness.
Its effect on the genus-one $N=5$ comparison is discussed explicitly in
\S\ref{subsec:diag_validation_extension}; in particular, disabling it removes the remaining convention-dependent discrepancy
with the published skeleton list.
\end{remark}

\begin{remark}[Canonicalization conventions used for comparison with published tables]
\label{rem:compare_canon_rules}
The invariants $\langle D\rangle(a,x)$ and $X_D(a,x)$ are defined intrinsically,
but published tables may incorporate harmless normalizations and equivalences.
When we compare our values to the genus-one tables
(e.g.\ \cite{AkimovaMatveev2012,AkimovaMatveevTarkaev2018,AkimovaMatveevTarkaev2020}),
we apply the following explicit comparison rules.

\begin{itemize}
\item \emph{Mirror identification.}
When a published table identifies mirror-related diagrams, we allow
\[
Q_D(a,x)\sim Q_D(a^{-1},x),
\]
which corresponds to the global crossing switch of Remark~\ref{rem:global_switch}.
At the skeleton level, this reverses each coefficient tuple $S_m(Q_D)$.

\item \emph{Overall $\pm a^k$ shifts (only when explicitly required).}
Some tabulations present polynomials after an overall normalization.
When necessary, we therefore allow multiplication by a global factor $\pm a^k$ and state explicitly when this convention is applied.
Such a shift does not affect the skeleton, since the $a$--exponents are discarded there.

\item \emph{Links.}
For links we avoid orientation dependence by using the bracket
$\langle D\rangle(a,x)$ as the primary key.
Alternatively, one may canonicalize $X_D$ over all component orientations; there are $2^{c(P)}$ such choices, but only
$2^{c(P)-1}$ distinct values up to simultaneous reversal of all component orientations.
We do not use this option in the $N=4,5$ genus-one validation.
\end{itemize}
\end{remark}

\subsection{Validation and extension}
\label{subsec:diag_validation_extension}

We conclude by recording validation checks against published genus-one tabulation data.
For the benchmark comparisons at $N\le 5$, the counts are the \emph{raw locally reduced} diagram classes
produced by the pipeline at crossing number $N$, that is, the within-level classes obtained from the
projection enumeration and the local diagram-level filters of \S\ref{subsec:diag_assignments_and_rules},
before removing classes already realized at smaller crossing numbers.
Optional comparison conventions from \S\ref{subsec:diag_skeleton_conventions} are mentioned separately when they affect the comparison.
The remaining discrepancies between our raw counts and the published tables are analyzed case by case below, and
their main sources are identified explicitly.

The purpose is to verify that the projection pool, diagram-level conventions, and state-sum
invariants are implemented correctly within the permutation model, with all comparison
conventions stated explicitly.

Table~\ref{tab:g1_diagram_counts} summarizes the comparison between our diagram enumeration
and the published genus-one tabulation data used as benchmarks (with sources indicated by superscripts).

\begin{table}[t]
\centering
\begin{tabular}{c|c|c}
\toprule
$N$ & knots & links \\
\midrule
$2$ & $1/1^{*}$ & $1/1^{*}$ \\
$3$ & $3/3^{*}$ & $4/4^{*}$ \\
$4$ & $18/17^{**}$ & $22/21^{***}$ \\
$5$ & $71/69^{**}$ & $99/99^{***}$ \\
\bottomrule
\end{tabular}
\caption{Comparison of \emph{raw locally reduced} diagram-class counts on the torus $T^2$ for $N\le 5$.
Each entry is written as $\text{this work}/\text{published count}$, where the left entry denotes the number of
within-level classes produced by our pipeline before ``new-at-$N$'' filtering.
Superscripts indicate the reference used for the published numbers:
${}^{*}$~\cite{Grishanov2009b} (early genus-one lists at $N\le 3$);
${}^{**}$~\cite{AkimovaMatveev2014};
${}^{***}$~\cite{AkimovaMatveevTarkaev2018,AkimovaMatveevTarkaev2020}.}
\label{tab:g1_diagram_counts}
\end{table}

\paragraph{The case $N=4$.}
For knots, our raw locally reduced count is $18$, whereas \cite{AkimovaMatveev2014} lists $17$ prime knots.
The extra class comes from the unique $N=4$ projection without bigon faces.
According to the published classification, this case is excluded by a global primeness convention on the torus that is not
captured by the local filters used here.
Thus the discrepancy reflects a genuinely global convention, rather than a failure of the local reduction stage.

For links, our raw count is $22$.
The published count $21$ is obtained from an original list of $22$ by a global equivalence already noted in
\cite{AkimovaMatveevTarkaev2020}.
Once that identified equivalence is taken into account, the comparison is consistent.

\paragraph{The case $N=5$.}
For links, our raw locally reduced count $99$ coincides exactly with \cite{AkimovaMatveevTarkaev2020}.
This comparison is made without imposing the optional over/under participation filter of
Definition~\ref{def:overunder_filter}; its role in the $N=5$ link comparison is discussed in
Remark~\ref{rem:overunder_convention}.

For knots at $N=5$, our raw locally reduced count is $71$, whereas \cite{AkimovaMatveev2014} lists $69$.
Of the two extra raw classes, one is equivalent, under the comparison conventions of
Remark~\ref{rem:compare_canon_rules}, to a knot type already present at smaller crossing number.
Thus enforcing the ``new-at-$N$'' convention lowers our count from $71$ to $70$.
The remaining discrepancy from $70$ to $69$ persists under the local diagram-level filters implemented in this paper.
Reconciling this final case with the published $N=5$ knot tables therefore requires an additional global convention beyond the
$\Omega_2$-based bigon rule, of the kind described in \cite{AkimovaMatveev2014}, such as the $\Omega_3$-type
central-fragment convention.

\begin{remark}[Reproducibility]
For each diagram in the computed lists we provide the unsensed canonical representative
of the underlying projection, the crossing assignment, and the corresponding state-sum invariants.
All counts and comparisons in this section are therefore directly checkable.
\end{remark}

\paragraph{Extension to $N=6,7,8,9,10$.}
\label{subsec:extension}

After validating the implementation for $N\le 5$, we extend the computation to
crossing numbers $N=6,7,8,9,10$, using the same invariants and comparison conventions throughout.
To our knowledge, these are the first complete genus-one tabulations at these
crossing numbers.
Under the comparison conventions of Remark~\ref{rem:compare_canon_rules}, the resulting \emph{new-at-$N$}
counts of knot and link types (knots/links) --- that is, the classes at crossing number $N$ that are not
equivalent to any class occurring at smaller crossing numbers --- are $378/525$ at $N=6$, $1743/2909$ at $N=7$,
$10704/16752$ at $N=8$, $61912/109231$ at $N=9$, and $434010/718823$ at $N=10$.
A detailed analysis of the data, including growth rates, structural patterns, and conjectures, is given in
Section~\ref{sec:observations}.

\section{Computational observations and conjectures}\label{sec:observations}

The tabulation data produced for $N\le 10$ contain considerably more information than
the bare counts reported in
Sections~\ref{sec:projections} and~\ref{sec:diagrams}.
In this section we examine the data systematically and record structural
patterns and natural conjectures arising from the computation.
Three results are proved
(Propositions~\ref{prop:lattice_lower_bound},
\ref{prop:max_bigons}, and~\ref{prop:span_parity}),
and several conjectures are formulated.
The empirical statements and examples supporting these conjectures can be
checked directly against the released datasets~\cite{OmelchenkoTorusMaps}.

Throughout, ``knot type'' and ``link type'' refer to equivalence classes of
diagrams under the canonical comparison conventions stated in
Remark~\ref{rem:compare_canon_rules}, and ``new at~$N$'' means not equivalent to any
diagram occurring at crossing numbers smaller than~$N$.

\subsection{Growth of genus-one diagram types}
\label{subsec:obs_growth}

Table~\ref{tab:obs_growth_summary} records the new-at-$N$ counts of
genus-one knot and link types and their ratio to the classical prime knot
count $K_0(N)$ (see~\cite{HosteThistlethwaiteWeeks1998}).

\begin{table}[ht]
\centering
\begin{tabular}{r|rrr|r|r}
\toprule
$N$ & knots & links & total & ratio to $N{-}1$ & knots$/K_0(N)$ \\
\midrule
$2$  & $1$      & $1$       & $2$        & ---   & --- \\
$3$  & $3$      & $4$       & $7$        & $3.5$ & $3$ \\
$4$  & $18$     & $22$      & $40$       & $5.7$ & $18$ \\
$5$  & $71$     & $99$      & $170$      & $4.2$ & $35.5$ \\
$6$  & $378$    & $525$     & $903$      & $5.3$ & $126$ \\
$7$  & $1743$   & $2909$    & $4652$     & $5.2$ & $249$ \\
$8$  & $10704$  & $16752$   & $27456$    & $5.9$ & $509.7$ \\
$9$  & $61912$  & $109231$  & $171143$   & $6.2$ & $1263.5$ \\
$10$ & $434010$ & $718823$  & $1152833$  & $6.7$ & $2630.4$ \\
\bottomrule
\end{tabular}
\caption{New genus-one diagram types at each crossing number~$N$.
The column ``ratio to $N{-}1$'' is the ratio of consecutive totals;
``knots$/K_0(N)$'' compares with the number of prime classical knots.}
\label{tab:obs_growth_summary}
\end{table}

Excluding the small-sample value at $N=3$, the consecutive ratios fluctuate
between $4.2$ and $6.7$.
A least-squares fit to $\log T_N$ over $N=3,\dots,10$, where
$T_N=\#\{\text{new types at }N\}$, gives
\begin{equation}\label{eq:obs_total_growth}
T_N \approx 0.0381\cdot 5.473^N
\qquad (R^2=0.9987 \text{ on the log scale}),
\end{equation}
suggesting an apparent exponential growth base between $5$ and $6$ in the computed range.
We do not make any asymptotic claim beyond this empirical fit.
Through $N=10$ the tabulation contains $508\,840$ knot and $848\,366$ link types
($1\,357\,206$ in total).
Compared with the classical prime-knot counts $K_0(N)$ in the same crossing range,
the genus-one knot counts are already larger by more than three orders of magnitude at
$N=10$, reflecting the combinatorial explosion that accompanies the passage from
genus~$0$ to genus~$1$.

\subsection{Maximum number of straight-ahead components}
\label{subsec:obs_max_components}

Recall from Section~\ref{sec:projections} that the number of straight-ahead components
$c(P)$ is computed from the permutation data by
Lemma~\ref{lem:straight_ahead_components}.
For each~$N$ let
\[
c_{\max}(N):=\max\{c(P): P \text{ is a prime unsensed projection with }N\text{ crossings}\}.
\]
Table~\ref{tab:obs_max_components} records the observed values.

\begin{table}[ht]
\centering
\begin{tabular}{r|c|c|c|l}
\toprule
$N$ & $c_{\max}(N)$ & $\lfloor N/2\rfloor+2$ & \# maximizers & face degrees of maximizer(s) \\
\midrule
$3$  & $3$ & $3$ & $1$ & $(3,3,6)$ \\
$4$  & $4$ & $4$ & $1$ & $(4,4,4,4)$ \\
$5$  & $4$ & $4$ & $1$ & $(3,3,4,4,6)$ \\
$6$  & $5$ & $5$ & $1$ & $(4,4,4,4,4,4)$ \\
$7$  & $5$ & $5$ & $1$ & $(3,3,4,4,4,4,6)$ \\
$8$  & $6$ & $6$ & $1$ & $(4,4,4,4,4,4,4,4)$ \\
$9$  & $6$ & $6$ & $2$ & $(3,3,4,4,4,4,4,4,6)$; $(4,4,4,4,4,4,4,4,4)$ \\
$10$ & $7$ & $7$ & $1$ & $(4,4,4,4,4,4,4,4,4,4)$ \\
\bottomrule
\end{tabular}
\caption{Maximum number of straight-ahead components among prime torus projections.
The formula $c_{\max}(N)=\lfloor N/2\rfloor+2$ holds throughout the computed range $N\le 10$.
The maximizer is unique except at $N=9$, where two maximizers occur.}
\label{tab:obs_max_components}
\end{table}

We first establish the lower bound by an explicit lattice construction.

\begin{proposition}[Lattice construction achieving $\lfloor N/2\rfloor+2$ components]
\label{prop:lattice_lower_bound}
For every even integer $N\ge 4$, let $m=N/2$ and let $P_N$ be the $4$--regular map
on the torus obtained from the square lattice $\mathbb{Z}^2$ by identifying under the
sublattice $\Lambda_N=\langle(2,0),(0,m)\rangle$.
Then:
\begin{enumerate}
\item $P_N$ is a connected equi-quadrilateral map with $N$ vertices and $N$ faces of
      degree~$4$;
\item the number of straight-ahead components is $c(P_N)=m+2=N/2+2$;
\item $P_N$ is prime in the sense of Definition~\ref{def:prime_projection_computational}:
      its underlying multigraph has no $2$--edge--cut, and $P_N$ is not split-witnessed.
\end{enumerate}
In particular, $c_{\max}(N)\ge \lfloor N/2\rfloor+2$ for all even $N\ge 4$.
\end{proposition}

\begin{proof}
The vertex set is $\mathbb{Z}^2/\Lambda_N=\{(i,j):0\le i<2,\;0\le j<m\}$,
so $|V|=2m=N$.
Edges connect lattice-adjacent vertices modulo $\Lambda_N$, giving $|E|=2N$;
since every face is a unit square, $|F|=N$ and $|V|-|E|+|F|=0$, confirming genus~$1$.

At each vertex the cyclic dart order is (right, up, left, down), so $\sigma^2$ pairs
right$\leftrightarrow$left and up$\leftrightarrow$down.
A straight-ahead walk alternates edge-traversal ($\alpha$) with the opposite-dart step
($\sigma^2$) and therefore traces a straight lattice line.
Horizontal lines $\{(\cdot,j):j\text{ fixed}\}$ close after $2$ steps
(since $(2,0)\in\Lambda_N$), producing $m$ horizontal components.
Vertical lines $\{(i,\cdot):i\text{ fixed}\}$ close after $m$ steps
(since $(0,m)\in\Lambda_N$), producing $2$ vertical components.
Hence $c(P_N)=m+2$.

To verify primeness, write the vertices as
\[
A_j=(0,j),\qquad B_j=(1,j),\qquad j\in \mathbb{Z}/m\mathbb{Z}.
\]
Every vertex is mixed: one horizontal and one vertical straight-ahead component meet there.
Each horizontal component meets exactly two mixed vertices, while each vertical component
meets $m$ mixed vertices. Hence $m(K)\ge 2$ for every straight-ahead component $K$, so
$P_N$ is not split-witnessed.

It remains to exclude $2$--edge--cuts in the underlying multigraph
$G(P_N)\cong C_2\square C_m$.
Let $S$ be a nonempty proper subset of vertices, and define
\[
I_A:=\{j: A_j\in S\},\qquad I_B:=\{j: B_j\in S\}\subset \mathbb{Z}/m\mathbb{Z}.
\]
Then the size of the edge cut $\partial S$ is
\[
|\partial S|=\beta(I_A)+\beta(I_B)+2|I_A\triangle I_B|,
\]
where $\beta(I)=0$ if $I=\varnothing$ or $I=\mathbb{Z}/m\mathbb{Z}$, and
$\beta(I)\ge 2$ otherwise, because the boundary of a nontrivial subset of a cycle
contains at least two edges.

If $I_A=I_B$, then $S$ is a nontrivial union of whole levels, so
$\beta(I_A)+\beta(I_B)\ge 4$.
If $I_A\neq I_B$, then either $|I_A\triangle I_B|\ge 2$, giving
$|\partial S|\ge 4$, or one of $I_A,I_B$ is a nontrivial proper subset, in which case
again $|\partial S|\ge 2+2=4$.
Thus every edge cut has size at least $4$, so $G(P_N)$ has no $2$--edge--cut.
Therefore $P_N$ is prime.
\end{proof}

\begin{conjecture}[Maximum component count]\label{conj:max_components}
For all $N\ge 3$, the maximum number of straight-ahead components among prime
$4$--regular torus projections with $N$ crossings satisfies
\begin{equation}\label{eq:obs_max_comp_formula}
c_{\max}(N)=\lfloor N/2\rfloor+2.
\end{equation}
\end{conjecture}

\begin{remark}[Status and observed structure of the maximizers]
Proposition~\ref{prop:lattice_lower_bound} establishes the lower bound
$c_{\max}(N)\ge\lfloor N/2\rfloor+2$ for even~$N$; for odd~$N$ the equality is verified
computationally for $N\le 9$.
The matching upper bound remains open in general.  In the computed range through $N=10$, the maximizer is unique
except at $N=9$.  For even~$N\le 10$, one maximizer is the lattice projection $P_N$ of
Proposition~\ref{prop:lattice_lower_bound}, whose faces are all quadrilaterals.  For odd
$N\le 9$, maximizers include a triangle--hexagon defect pattern with sorted face-degree vector
$(3,3,\underbrace{4,\ldots,4}_{N-3},6)$; at $N=9$ an additional equi-quadrilateral maximizer also appears.
An explicit construction and classification of the odd maximizers would complete the picture.
\end{remark}

\begin{remark}[The maximizer count $\lfloor N/2\rfloor+2$ versus naive bounds]
A trivial upper bound is $c(P)\le 2N$, since $c(P)=\tfrac12\,\#\mathrm{cyc}(\rho)$
and $\rho$ acts on $4N$ darts.
The observed bound $\lfloor N/2\rfloor+2$ is dramatically smaller, indicating that
the torus topology and the primeness constraint impose strong restrictions on how many
independent components can coexist in a single projection.
\end{remark}

\subsection{Equi-quadrilateral projections are always links}
\label{subsec:obs_equiquad}

A $4$--regular map on the torus is called \emph{equi-quadrilateral} if every face has
degree~$4$.  In the notation of Section~\ref{sec:maps}, this means that every cycle of
the face permutation $\varphi=\sigma\alpha$ has length exactly~$4$.

\begin{table}[ht]
\centering
\begin{tabular}{r|r|l}
\toprule
$N$ & \# equi-quad. & component counts \\
\midrule
$3$  & $1$ & $2$ \\
$4$  & $3$ & $2, 3, 4$ \\
$5$  & $2$ & $2, 2$ \\
$6$  & $4$ & $2, 3, 4, 5$ \\
$7$  & $2$ & $2, 2$ \\
$8$  & $6$ & $2, 2, 3, 4, 5, 6$ \\
$9$  & $4$ & $2, 2, 4, 6$ \\
$10$ & $6$ & $2, 2, 3, 3, 6, 7$ \\
\bottomrule
\end{tabular}
\caption{Equi-quadrilateral prime projections and their straight-ahead component counts.
In every case $c(P)\ge 2$.}
\label{tab:obs_equiquad}
\end{table}

\begin{conjecture}[No equi-quadrilateral knot projections]\label{conj:equiquad_links}
Every prime equi-quadrilateral $4$--regular map on the torus
has at least $2$ straight-ahead components.  In particular, no equi-quadrilateral
projection supports a knot diagram.
\end{conjecture}

\begin{remark}[Algebraic reformulation and partial evidence]
\label{rem:equiquad_reformulation}
The equi-quadrilateral condition states that the face permutation
$\varphi=\sigma\alpha$ is itself a product of $N$ disjoint $4$--cycles, so
$\sigma$ and $\varphi$ share the same cycle type.
Since $\rho=\sigma^2\alpha=\sigma\varphi$, the conjecture admits a purely
permutation-theoretic reformulation: if $\sigma$ and $\varphi$ are
permutations of $\{1,\dots,4N\}$, both of cycle type $(4^N)$, such that
\[
\alpha:=\sigma^{-1}\varphi
\]
is a fixed-point-free involution and
$\langle\sigma,\alpha\rangle$ (equivalently, $\langle\sigma,\varphi\rangle$)
acts transitively, must the product $\rho=\sigma\varphi$ have at least $4$ cycles?
Equivalently, can $\rho$ have as few as $2$ cycles?

The component-count evidence in Table~\ref{tab:obs_equiquad} verifies Conjecture~\ref{conj:equiquad_links} through $N\le 10$.
In the range $N\le 8$, a stronger structural pattern was also checked:
every straight-ahead component of every equi-quadrilateral projection represents
a nonzero class in $H_1(T^2;\mathbb{Z})\cong\mathbb{Z}^2$.
The Euler--Hurwitz formula for the factorization $\rho=\sigma\varphi$ gives
\[
  \#\mathrm{cyc}(\rho)=2N+2-2g',
\]
where $g'\ge 0$ is the genus of the covering surface associated with the
transitive action of $\langle\sigma,\varphi\rangle$ on $\{1,\dots,4N\}$
in the sense of the Riemann--Hurwitz formula for permutation factorisations
(see, e.g., \cite{LandoZvonkin2004}).
A single-component projection ($c(P)=1$, i.e.\ $\#\mathrm{cyc}(\rho)=2$)
would require $g'=N$, the maximum possible value.
Showing that the torus map structure (connectivity, cellularity, and the involution
constraint on~$\alpha$) forces $g'<N$ would settle the conjecture.
\end{remark}

\subsection{Coefficient structure of the genus-one bracket}
\label{subsec:obs_bracket_structure}

\paragraph{The $a$-span.}
For a knot diagram $D$ on the torus, define the \emph{$a$-span} of the normalized
polynomial $X_D(a,x)$ as the difference between the maximum and minimum exponents
of~$a$ appearing with nonzero coefficient:
\[
\operatorname{span}_a(X_D):=\max\{k: [a^k]X_D\neq 0\}
  -\min\{k: [a^k]X_D\neq 0\}.
\]
Table~\ref{tab:obs_span} records the observed extremes and averages.

\begin{table}[ht]
\centering
\begin{tabular}{r|c|c|c|c}
\toprule
$N$ & min span & max span & $4N$ & average span \\
\midrule
$3$  &  $6$ & $12$ & $12$ & $9.0$ \\
$4$  &  $8$ & $16$ & $16$ & $11.9$ \\
$5$  & $10$ & $20$ & $20$ & $15.4$ \\
$6$  &  $4$ & $24$ & $24$ & $17.7$ \\
$7$  &  $4$ & $28$ & $28$ & $21.0$ \\
$8$  &  $4$ & $32$ & $32$ & $23.6$ \\
$9$  &  $4$ & $36$ & $36$ & $26.7$ \\
$10$ &  $4$ & $40$ & $40$ & $29.3$ \\
\bottomrule
\end{tabular}
\caption{Observed $a$-span of the normalized polynomial $X_D(a,x)$ for genus-one knot types.
The maximum span equals $4N$ at every crossing number in the computed range $3\le N\le 10$.}
\label{tab:obs_span}
\end{table}

\begin{proposition}[Parity of the $a$-span]\label{prop:span_parity}
Let $D=(P,b)$ be a diagram on the torus with $N$ crossings (knot or link).
Then every monomial in the bracket $\langle D\rangle(a,x)$ has $a$-exponent of the
same parity as~$N$.
In particular, $\operatorname{span}_a\langle D\rangle$ is an even integer.
For knots, $X_D(a,x)=(-a)^{-3w(D)}\langle D\rangle(a,x)$ shifts all
$a$-exponents by a constant, so $\operatorname{span}_a(X_D)=
\operatorname{span}_a\langle D\rangle$ is again even.
\end{proposition}

\begin{proof}
In the state sum
(Definition~\ref{def:generalized_bracket_section})
\[
\langle D\rangle(a,x)=\sum_{s:\mathrm{Vert}(P)\to\{A,B\}}
a^{n_A(s)-n_B(s)}\bigl(-a^2-a^{-2}\bigr)^{\gamma(s)}x^{\delta(s)},
\]
the weight exponent satisfies $n_A(s)-n_B(s)=2n_A(s)-N\equiv N\pmod{2}$
for every state~$s$, since $n_A(s)+n_B(s)=N$.
Moreover, $(-a^2-a^{-2})^{\gamma(s)}$ is a Laurent polynomial in~$a$
whose every term has even $a$-exponent.
Therefore the $a$-exponent of each monomial in $\langle D\rangle$ has
parity~$N$; in particular, the maximum and minimum $a$-exponents share the
same parity, so $\operatorname{span}_a\langle D\rangle$ is even.

For knots, $X_D(a,x)=(-a)^{-3w(D)}\langle D\rangle(a,x)$ shifts all
$a$-exponents by the constant~$-3w(D)$, leaving the span unchanged.
\end{proof}

\begin{conjecture}[$a$-span bound]\label{conj:span}
For every genus-one knot diagram $D$ with $N$ crossings,
\[
\operatorname{span}_a(X_D)\le 4N,
\]
and this bound is achieved at every~$N$ in the computed range $3\le N\le 10$.
The same upper bound holds for the bracket $\langle D\rangle(a,x)$ of link
diagrams for $3\le N\le 10$.
\end{conjecture}

\begin{remark}[Analogy with the classical span theorem]
In the classical setting ($\Sigma_0=S^2$), the Kauffman--Murasugi--Thistlethwaite
theorem states that the $a$-span of the Jones polynomial of a reduced alternating link
diagram with $N$ crossings equals exactly~$4N$ (see \cite{Kauffman1999} for the
virtual-knot perspective).
Conjecture~\ref{conj:span} asserts that the same numerical bound $4N$ serves as an
upper bound in the genus-one case, though the underlying diagrams are typically not
alternating in the surface-based sense.
In the computed range the average $a$-span grows steadily, while the maximum remains exactly $4N$.
For $N\ge 6$ the detailed $a$-span distribution develops a bimodal character:
a primary peak near $\operatorname{span}_a\approx 2.5N$--$3N$ and a secondary
concentration at the maximum $\operatorname{span}_a=4N$.
An extension of the classical span inequality to the generalized bracket on surfaces
of positive genus would provide a theoretical foundation for this observed bound.
\end{remark}

\paragraph{Bound on essential state circles.}
The variable~$x$ likewise obeys a tight bound in the analyzed range.

\begin{conjecture}[Bound on essential state circles]\label{conj:xmax}
For every genus-one knot diagram $D$ with $N$ crossings, the maximum power
of~$x$ in $X_D(a,x)$ satisfies $x_{\max}(D)\le N-1$.
This bound is achieved at every~$N$ in the computed range $3\le N\le 10$.
\end{conjecture}

\paragraph{Purely essential knots.}

\begin{definition}[Purely essential knot]\label{def:purely_essential}
A genus-one knot type is called \emph{purely essential} if its invariant
$X_D(a,x)$ has no $x^0$~term.
This is a genuinely genus-one phenomenon, vacuous on~$S^2$.

A sufficient diagrammatic condition for this is that every smoothing state
contains at least one essential (non-contractible) circle; the converse need
not hold a priori, since the $x^0$-contributions of different states may cancel.
\end{definition}

\begin{table}[ht]
\centering
\begin{tabular}{r|rr|r}
\toprule
$N$ & purely essential & total & fraction \\
\midrule
$3$  &     $2$ &      $3$ & $0.67$ \\
$4$  &    $14$ &     $18$ & $0.78$ \\
$5$  &    $34$ &     $71$ & $0.48$ \\
$6$  &   $276$ &    $378$ & $0.73$ \\
$7$  &  $1053$ &   $1743$ & $0.60$ \\
$8$  &  $7805$ &  $10704$ & $0.73$ \\
$9$  & $42789$ &  $61912$ & $0.69$ \\
$10$ & $325551$ & $434010$ & $0.75$ \\
\bottomrule
\end{tabular}
\caption{Fraction of purely essential knot types at each crossing number.
For even~$N$ the fraction is consistently around~$0.73$--$0.78$ in the computed range.}
\label{tab:obs_purely_essential}
\end{table}

\begin{observation}[Dominance of purely essential knots]\label{obs:purely_essential}
For $N\ge 4$, a substantial fraction --- and for even~$N$ the majority --- of genus-one
knot types are purely essential. The fraction shows a pronounced parity dependence:
for even~$N\in\{4,6,8,10\}$ it stays in the range $0.73$--$0.78$, while for odd
$N\in\{5,7,9\}$ it ranges from $0.48$ to $0.69$.
\end{observation}

\subsection{Bigon prevalence}
\label{subsec:obs_bigons}

\begin{proposition}[Maximum number of bigon faces]\label{prop:max_bigons}
A $4$--regular map on the torus with $N$ vertices has at most $N-1$ bigon faces.
When this bound is attained, the unique non-bigon face has degree~$2N+2$.
\end{proposition}

\begin{proof}
By the Euler formula, $F=N$.
The sum of face degrees equals the number of darts, $4N$.
If all $N$ faces were bigons, this sum would be $2N<4N$, a contradiction;
hence at most $N-1$ faces have degree~$2$.
If exactly $N-1$ are bigons, the remaining face has
degree $4N-2(N-1)=2N+2$.
\end{proof}

\begin{table}[ht]
\centering
\begin{tabular}{r|rr|r|r}
\toprule
$N$ & with bigons & prime total & fraction & max \# bigons \\
\midrule
$3$  &     $3$ &      $6$ & $0.50$ & $2$ \\
$4$  &    $17$ &     $23$ & $0.74$ & $3$ \\
$5$  &    $68$ &     $81$ & $0.84$ & $4$ \\
$6$  &   $342$ &    $379$ & $0.90$ & $5$ \\
$7$  &  $1708$ &   $1795$ & $0.95$ & $6$ \\
$8$  &  $9536$ &   $9815$ & $0.97$ & $7$ \\
$9$  & $54749$ &  $55593$ & $0.98$ & $8$ \\
$10$ & $328325$ & $331233$ & $0.99$ & $9$ \\
\bottomrule
\end{tabular}
\caption{Bigon prevalence among prime projections.  The last column
confirms sharpness of the bound $N-1$ from
Proposition~\ref{prop:max_bigons} at every~$N$ in the computed range.}
\label{tab:obs_bigons}
\end{table}

\begin{observation}[Near-ubiquity of bigons]\label{obs:bigons}
The fraction of prime projections containing at least one bigon face increases
monotonically with~$N$, reaching $99\%$ at $N=10$.
The bound of Proposition~\ref{prop:max_bigons} is sharp at every $N\le 10$:
the maximum $N-1$ is achieved by exactly one prime projection
(with face-degree vector $(2^{N-1},2N{+}2)$).
\end{observation}

\paragraph{Closing remarks.}
Among the patterns suggested by the computation through $N=10$, the most robust are the component-count
conjecture \eqref{eq:obs_max_comp_formula}, the absence of equi-quadrilateral knot
projections (Conjecture~\ref{conj:equiquad_links}), and the span bounds of
Conjectures~\ref{conj:span} and~\ref{conj:xmax}.
At the projection level, Proposition~\ref{prop:max_bigons} and the observed near-ubiquity
of bigons suggest further asymptotic questions.
All empirical data underlying these observations are available in machine-readable form
in the released datasets~\cite{OmelchenkoTorusMaps}.

\section{Outlook beyond the torus}\label{sec:outlook}

The map-theoretic encoding of projections by permutation pairs $(\alpha,\sigma)$ is
not specific to genus one: it applies without modification on any closed orientable
surface $\Sigma_g$ and yields a uniform framework for generating cellular
$4$--regular projections up to (un)sensed equivalence.  The torus case, however,
is the first positive-genus setting in which the subsequent diagram stage remains
both homologically meaningful and computationally controllable.

\medskip\noindent
\emph{Why the torus case is tractable.}
On $T^2$, every embedded essential circle represents a primitive class in
$H_1(T^2;\mathbb Z)\cong\mathbb Z^2$ and hence has nonzero reduction in
$H_1(T^2;\mathbb F_2)$ (cf.\ Lemma~\ref{lem:cut_via_face_span_section}).
Consequently, the state-circle dichotomy \emph{contractible vs.\ essential} can be
detected by a simple $\mathbb F_2$ boundary test, and the standard genus-one
bracket used in the literature depends only on whether a state circle is trivial
or not, so that all nontrivial homology classes can be counted by a single
variable $x$.  This also makes the ambiguity coming from the mapping class group
action on $H_1(T^2)$ irrelevant after the standard genus-one specialization.

\medskip\noindent
\emph{What changes for $g\ge 2$.}
Separating state circles can be essential but null-homologous, so any
homology-refined state sum must keep track of integral classes in
$H_1(\Sigma_g;\mathbb Z)$, not merely ``zero vs.\ nonzero''.  Moreover,
comparison of such homology-labelled polynomials across diagrams naturally
involves the action of the mapping class group of $\Sigma_g$ on $H_1(\Sigma_g)$,
which factors through the surjection
${\rm MCG}(\Sigma_g)\twoheadrightarrow{\rm Sp}(2g,\mathbb Z)$.
Developing efficient canonicalization procedures for this additional symmetry,
or working with mapping-class-group--invariant specializations of higher-genus
state-sum invariants, is an important direction that we leave for future work.

\medskip\noindent
\emph{Computational outlook.}
Cellularity forces $N\ge 2g-1$ (since $F=N+2-2g\ge 1$), and the state sum remains
exponential in~$N$; extending the computation to higher genus will therefore
require both algebraic and algorithmic innovations beyond those developed here.
The genus-one case treated in this paper, with its complete tabulation through
$N=10$ and the structural patterns it has revealed, provides both the baseline
data and the methodological template for such extensions.

\paragraph{Data availability.}
The reference implementation and all machine-readable datasets underlying this
study --- projection lists, diagram lists, and invariant tables through
$N=10$ --- are openly available in the public repository
\texttt{Torus\_Maps} at \url{https://github.com/avo-travel/Torus_Maps}~\cite{OmelchenkoTorusMaps}.
Every table, count, and observation in the paper can be reproduced from the
released code and data.

\appendix

\section{Algorithmic realization of projection enumeration}\label{app:proj_algorithms}

This appendix records the computational realization of the projection enumeration step of
Section~\ref{sec:projections}.  The permutation model itself is the classical map formalism recalled in
Section~\ref{sec:maps}.  The points described here are implementation choices: they make the enumeration
fully reproducible and explain why the fixed-torus projection search is feasible in practice.

The general orbit-enumeration viewpoint for permutation-encoded maps is classical; in the closely related
setting of immersed circles, Coquereaux--Zuber describe the problem as an orbit enumeration under centralizer
actions and use both direct orbit construction and more economical group-theoretic reductions
\cite{CoquereauxZuber2016}.  Our implementation is adapted to a different task: the construction of a
canonical pool of prime knot and link projections in the fixed torus.  We therefore keep the vertex rotation
fixed and enumerate only the edge-pairing involution, while using several safe look-ahead tests during the
backtracking.

\subsection{Data model and recursive state}
We represent a labelled map by permutations $(\alpha,\sigma)$ on
$H=\{1,\dots,4N\}$, stored as arrays of images.  Throughout the enumeration we fix
$\sigma=\sigma_0$ as in \eqref{eq:standard_sigma_proj} and generate $\alpha$ as a
fixed-point-free involution, equivalently as a perfect matching of $H$.

The recursive state consists of the partial involution $\alpha$, the set of unused darts, the set of activated
vertices, and two pieces of topological look-ahead data.  First, a rollback disjoint-set structure is maintained
on the $N$ vertex blocks, so that the number $q$ of connected components of the partial abstract graph is known
at every node.  Second, the implementation maintains the current boundary permutation of the partial ribbon
surface obtained after attaching the already chosen edge bands.  Initially this boundary permutation is
$\sigma_0$; if the new pair is $\alpha(i)=j$, then the boundary permutation is updated by swapping the two
successors of $i$ and $j$.  Thus the update is local and reversible.

For speed, dart sets are stored as bit masks: bit $d-1$ records the presence of the dart $d$.  Since the largest
values used in the tables have $4N\le 40$ darts, the active and unused dart sets fit into a single machine word.
This is only an implementation device; it does not change the mathematical search tree.

\subsection{Canonical construction path for matchings}
Naively, the number of perfect matchings of $H$ is $(4N-1)!!$.  We reduce the main relabelling redundancy by
using the following construction rule.

\smallskip
\noindent\textbf{Rule (activated vertices).}
At each stage of the recursion there is an initial segment of activated vertices in the block order induced by
$\sigma_0$.  Let $i$ be the smallest unused dart.  We may pair $i$ either
(a) with an unused dart on an activated vertex, or
(b) with the first dart of the first non-activated vertex, thereby activating that vertex.
This fixes the first appearance of every new vertex.  Consequently the search contains at least one labelled
representative of each unlabelled map with $\sigma=\sigma_0$, while avoiding many relabelled copies.

The loop-edge and monogon restrictions are enforced immediately.  If loop edges are forbidden, candidates in the
same vertex block as $i$ are removed.  If monogons are forbidden, the adjacent pairs
$\{i,\sigma_0(i)\}$ and $\{i,\sigma_0^{-1}(i)\}$ are removed; these are exactly the local pairings that create a
fixed point of the face permutation $\varphi=\sigma_0\alpha$.

\subsection{Dynamic boundary and genus look-ahead}
Let $m$ be the number of already chosen $\alpha$-pairs and let
\[
        r=2N-m
\]
be the number of still missing edges.  Let $B$ be the number of cycles of the current boundary permutation and
let $q$ be the number of connected components of the partial abstract graph on the $N$ vertex blocks.
The target number of faces for genus $g$ is
\[
        F_g=N+2-2g.
\]
For the torus, $g=1$ and $F_g=N$.

When a new pair $\alpha(i)=j$ is added, the two darts $i,j$ either lie on the same current boundary cycle or on
two distinct boundary cycles.  In the first case the attachment splits one boundary component and
$B$ increases by $1$; in the second case it merges two boundary components and $B$ decreases by $1$.
Thus each remaining edge can change $B$ only by $\pm1$.  This gives the face-count and parity tests
\begin{equation}\label{eq:app_face_pruning}
        |B-F_g|\le r,
        \qquad
        B-F_g\equiv r \pmod 2.
\end{equation}
The rollback disjoint-set structure gives the connectedness look-ahead
\begin{equation}\label{eq:app_conn_pruning}
        q-1\le r,
\end{equation}
since each future edge can reduce the number of connected components by at most one.

Finally, the partial ribbon surface has Euler characteristic
\[
        N-m = 2q-2g_{\mathrm{part}}-B,
\]
and hence
\begin{equation}\label{eq:app_partial_genus}
        g_{\mathrm{part}}=\frac{2q-B-N+m}{2}.
\end{equation}
A branch is rejected if $g_{\mathrm{part}}>g$.  Conversely, after spending at least $q-1$ future edges to connect
the remaining components, at most $r-(q-1)$ future edges remain available to increase genus.  Therefore we also
reject if
\begin{equation}\label{eq:app_genus_capacity}
        g_{\mathrm{part}}+r-(q-1)<g.
\end{equation}
All tests in \eqref{eq:app_face_pruning}--\eqref{eq:app_genus_capacity} are necessary conditions only.  Hence they
can prune the search tree without deleting any valid torus projection.

\subsection{Unsensed canonical representative}
Definition~\ref{def:proj_canonical_rep} defines the canonical representative as the minimum over an $S_H$-orbit.
The implementation computes the same invariant by rooted normalization.  Given a connected encoding
$(\alpha,\sigma)$ and a root dart $r\in H$, let $\mathrm{Norm}_r(\alpha,\sigma)$ denote the relabelled encoding
obtained by the fixed traversal routine used in the implementation: starting from $r$, inspect neighbors in the
order $(\sigma(h),\alpha(h))$ at each visited dart $h$, assign labels $1,2,\dots,|H|$ in order of first discovery,
and conjugate $(\alpha,\sigma)$ by the resulting relabelling permutation.  Thus
$\mathrm{Norm}_r(\alpha,\sigma)=(\alpha^{(r)},\sigma^{(r)})$.

The unsensed canonical representative is obtained by taking the minimum over all roots and over the two surface
orientations:
\[
\mathrm{Can}_{\mathrm{unsensed}}(\alpha,\sigma)
=\min\Bigl(
\{\mathrm{Norm}_r(\alpha,\sigma)\}_{r\in H}
\cup
\{\mathrm{Norm}_r(\alpha,\sigma^{-1})\}_{r\in H}
\Bigr),
\]
with respect to the fixed lexicographic order on the pair $(\alpha^{(r)},\sigma^{(r)})$.

The optimized implementation preserves this canonical key.  It only avoids building full rooted normal forms that
cannot be minimal.  Since the first coordinate of the key is $\alpha^{(r)}(1)$, i.e. the label assigned to
$\alpha(r)$ in the rooted traversal, we first compute this single value for every root and both orientations.  Full
normalization is then performed only for roots attaining the minimum possible first coordinate.  This is a safe
prefilter: a root with larger $\alpha^{(r)}(1)$ cannot minimize the full lexicographic key.  The resulting canonical
representatives are therefore identical to those obtained by the full rooted search.

\subsection{Parallel frontier}
The enumeration is also parallelized in a deterministic way.  A short initial segment of the backtracking tree is
explored serially until a prescribed depth $k$ is reached.  Each resulting prefix is stored simply as a list of
chosen pairs $(i,j)$.  A worker process replays this prefix, reconstructs the partial involution, the boundary
permutation and the rollback disjoint-set state, and then completes the suffix search independently.  Each worker
returns a local set of unsensed canonical representatives.  The main process takes the union of these sets and sorts
it before applying the projection-level filters and writing the output files.

This parallelization changes neither the construction rule nor the canonical representative.  It is only a way of
partitioning the same finite search tree into independent subtrees.

\subsection{Primeness witnesses and output}
For each unsensed canonical representative we compute:
(i) straight-ahead components via $\rho=\sigma^2\alpha$ (Lemma~\ref{lem:straight_ahead_components});
(ii) the $2$--edge--cut witness on the abstract multigraph $G(P)$ (Definition~\ref{def:composite_witness});
(iii) the splitness witness based on mixed vertices (Definition~\ref{def:split_witness}).
Filtering by these witnesses yields the ``prime total'' counts in Table~\ref{tab:g1_projection_counts} and the component
distribution in Table~\ref{tab:g1_links_by_components}, in the sense of Definition~\ref{def:prime_projection_computational}.
The released implementation writes the resulting prime projections both as human-readable permutation passports and
as deterministic machine-readable JSON files~\cite{OmelchenkoTorusMaps}.

\section{Algorithmic realization of diagram enumeration and state-sum evaluation}
\label{app:diagram_algorithms}

This appendix summarizes the computational realization of the diagram-level enumeration
described in Section~\ref{sec:diagrams}.
As in Appendix~\ref{app:proj_algorithms}, the purpose is not to introduce new mathematics,
but to make the pipeline fully reproducible and to document the algorithmic decomposition
used in the reference implementation.

\subsection{Input data: prime projections}
The diagram enumeration stage takes as input the lists of unsensed prime projections
constructed in Section~\ref{sec:projections}.
Algorithmically, each projection is specified by its unsensed canonical permutation encoding
$(\alpha,\sigma)$ together with the decomposition into straight-ahead components
(the $\langle \alpha,\sigma^2\rangle$--orbits).

Crucially, projection enumeration and diagram enumeration are completely decoupled:
the projection lists are computed once and stored in machine-readable form.
All subsequent diagram-level computations operate exclusively on these stored encodings
and do not re-enumerate projections.
This separation is essential for computational scalability at $N\ge 6$.

\subsection{Crossing assignments and global symmetries}
For a fixed projection $P$ with $N$ vertices, a diagram is determined by a crossing
assignment $b:\mathrm{Vert}(P)\to\{0,1\}$.
Naively this yields $2^N$ diagrams supported by $P$.

When the optional global crossing switch is enabled (Remark~\ref{rem:global_switch};
cf.\ Remark~\ref{rem:compare_canon_rules}), i.e.\ when the chosen comparison convention identifies
$b$ and $1-b$, we quotient by the involution $b\mapsto 1-b$ by fixing one reference bit of the crossing
assignment.  In the implementation this is treated as one additional linear equation, for example $b(v_1)=0$.

\subsection{Local diagram-level filters as linear constraints}
Before evaluating invariants, crossing assignments are subjected to the explicit diagram-level filters described in
Section~\ref{sec:diagrams}.  The optimized implementation uses the fact that the bigon rule is linear over
$\mathbb F_2$.  Each two-vertex bigon contributes an equation
\[
        b(u)+b(v)=\varepsilon \pmod 2,
\]
where $\varepsilon$ is determined by the local parities in Lemma~\ref{lem:bigon_omega2}.  The global crossing-switch
normalization, when enabled, contributes the additional equation $b(v_1)=0$.

Thus the admissible crossing assignments form an affine subspace of $\mathbb F_2^N$.  For each projection, the code
solves this system once by Gaussian elimination over $\mathbb F_2$ and enumerates only the resulting affine solution
space.  If the system is inconsistent, the projection contributes no diagrams under the chosen local rules.  This
replaces the older procedure of scanning all $2^N$ assignments and testing the bigon equations one by one.

The over/under participation filter for links is applied next.  For link diagrams ($c(P)\ge 2$), the implementation
scans the mixed vertices and records, for each straight-ahead component, whether it appears both as an overpassing and
as an underpassing component at mixed crossings.  This test is global but purely combinatorial and is still much cheaper
than a state-sum evaluation.

\subsection{State-sum structure and precomputation}
The generalized bracket $\langle D\rangle(a,x)$ is defined as a sum over all states
$s:\mathrm{Vert}(P)\to\{A,B\}$.
For implementation we encode a state by bits, writing $s(v)=0$ for the $A$--smoothing and
$s(v)=1$ for the $B$--smoothing, so that states are bit vectors in $\{0,1\}^N$.
A direct evaluation would require $2^N$ state traversals \emph{for each} crossing
assignment $b$, which quickly becomes prohibitive for $N\ge 6$.

The implementation uses the following reindexing.  Throughout this appendix the projection is given in the standardized
encoding $(\alpha,\sigma_0)$ of \S\ref{subsec:proj_model_constraints}, so that each vertex block has a canonical first dart.
With this convention, write the local cyclic order at the $k$th vertex as
$(h_0\,h_1\,h_2\,h_3)=(4k-3,\,4k-2,\,4k-1,\,4k)$ and define
\[
\tau_v^{(0)}:=(h_0\,h_1)(h_2\,h_3), \qquad
\tau_v^{(1)}:=(h_1\,h_2)(h_3\,h_0).
\]
Under the bit encoding of Definition~\ref{def:diagram_on_projection}, if $b(v)\in\{0,1\}$ is the crossing bit and
$s(v)\in\{0,1\}$ is the state bit, then the local smoothing selected at $v$ is
$\tau_v^{(b(v)\oplus s(v))}$.  Thus the smoothing geometry depends only on
\[
        t(v):=b(v)\oplus s(v).
\]
Since, for fixed $b$, the map $s\mapsto t=b\oplus s$ is a bijection of $\{0,1\}^N$, we precompute the geometric data once
for all $t\in\{0,1\}^N$ and then reuse them for every admissible crossing assignment $b$.

Concretely, for a fixed projection $P$ the code precomputes for every $t\in\{0,1\}^N$ the pair
\[
        (\gamma(t),\delta(t)),
\]
where $\gamma(t)$ and $\delta(t)$ are the numbers of contractible and essential state circles in the smoothing defined by
$t$.  The optimized implementation does not explicitly build the smoothing involution $\tau_t$ or the full state
permutation $\pi_t$ for every state.  Instead it precomputes two local partner arrays, corresponding to
$\tau_v^{(0)}$ and $\tau_v^{(1)}$, and follows the map
\[
        p_t(d)=\tau_t(\alpha(d))
\]
on the fly.  During this traversal it simultaneously records the $\mathbb F_2$ edge-incidence vector of each physical
state circle and tests it against the face-boundary span of Lemma~\ref{lem:cut_via_face_span_section}.  This gives
$\gamma(t)$ and $\delta(t)$ without constructing intermediate permutation objects.

This precomputation requires $2^N$ smoothing traversals once per projection and is independent of the particular crossing
assignment $b$.

\subsection{Fast evaluation for varying crossing assignments}
Given the precomputed table $\{(\gamma(t),\delta(t))\}_{t\in\{0,1\}^N}$, the bracket for
a specific diagram $D=(P,b)$ is evaluated by a single pass over all $t\in\{0,1\}^N$:
\[
\langle D\rangle(a,x)
=\sum_{t\in\{0,1\}^N}
a^{N-2|b\oplus t|}(-a^2-a^{-2})^{\gamma(t)}x^{\delta(t)}.
\]
No additional traversal of the smoothed diagram is required at this stage.

Several constant-factor optimizations are used in this pass.  The exponent
$N-2|b\oplus t|$ is stored in a table indexed by the XOR pattern $b\oplus t$, so the exponent table has size
$O(2^N)$ rather than $O(4^N)$.  Powers of the loop factor $(-a^2-a^{-2})^\gamma$, shifted by the possible state
exponents, are also precomputed.  For knot projections, the writhe of every crossing assignment is precomputed once
per projection and then used as a table lookup in the normalization
$X_D(a,x)=(-a)^{-3w(D)}\langle D\rangle(a,x)$.

As a result, for each projection the total cost is $O(2^N)$ for smoothing precomputation plus $O(2^N)$ per admissible
crossing assignment, but with substantially smaller Python overhead and memory use than the direct implementation.

\subsection{Parallelization over projections}
The diagram stage is parallelized over the input projections.  The projection list is divided into chunks, and each
worker process computes the local knot and link keys for its assigned projections.  The main process then merges the
worker outputs deterministically: for each canonical key it keeps the earliest representative in the original projection
order and crossing-assignment order.  Thus the parallel run produces the same representatives as the serial run, while
changing only the order in which independent projection chunks are processed.

\subsection{Canonicalization and incremental enumeration}
Diagrams are classified by canonical keys derived from the relevant state-sum invariant
($X_D(a,x)$ for knots and $\langle D\rangle(a,x)$ for links), together with the chosen equivalence conventions
(e.g.\ mirror normalization, cf.\ Remark~\ref{rem:compare_canon_rules}).
To support incremental ``new-at-$N$'' enumeration, the implementation maintains
\begin{enumerate}
\item[(i)] a persistent library of canonical keys coming from all smaller crossing numbers, and
\item[(ii)] a temporary set of keys already accepted at the current crossing number~$N$.
\end{enumerate}
A diagram encountered at level $N$ is recorded as new if and only if its canonical key belongs to neither set;
once accepted, the key is inserted into the level-$N$ set and merged into the persistent library after the pass
at crossing number $N$ completes.

For the largest layer reported here, $N=10$, the implementation processed $331\,233$ prime projections.
After the global crossing-switch convention, this corresponds to $169\,591\,296$ crossing assignments tested.
The bigon and over/under participation filters reduced these to $26\,833\,542$ diagrams for which polynomial keys were evaluated.
The layer produced $508\,229$ distinct knot keys and $845\,636$ distinct link keys at crossing number $10$,
of which $434\,010$ knot keys and $718\,823$ link keys were new relative to all smaller crossing numbers.
On the machine used for the final run, this layer required $680.26$ seconds.

\subsection{Reference code and datasets}
The implementation and accompanying datasets are available at~\cite{OmelchenkoTorusMaps}.
The repository contains the scripts for both projection enumeration
(Appendix~\ref{app:proj_algorithms}) and diagram evaluation described above,
together with machine-readable projection, diagram, and invariant tables through $N=10$.

\newpage
\bibliographystyle{unsrt}

\end{document}